\documentclass[review]{elsarticle}
\usepackage{amsmath,amsthm,amssymb,mathtools,amscd}
\usepackage{latexsym,graphicx,color,indentfirst}

\usepackage[linesnumbered,boxed,norelsize]{algorithm2e}
\usepackage{fancyhdr}
\usepackage{mathrsfs}
\usepackage{algorithmic}
\usepackage{subcaption}
\usepackage{float}
\usepackage{multirow}
\usepackage{array,xcolor,colortbl}
\usepackage{ctable} 
\usepackage{hyperref,amsmath,amssymb,amscd}

\makeatletter \@addtoreset{equation}{section} \makeatother
\newtheorem{theorem}{Theorem}[section]
\newtheorem{lemma}[theorem]{Lemma}
\newtheorem{assumption}[theorem]{Assumption}

\newtheorem{proposition}[theorem]{Proposition}

\newtheorem{example}[theorem]{Example}

\renewcommand{\theequation}{\arabic{section}.\arabic{equation}}
\renewcommand{\theenumi}{(\roman{enumi})}

\newcommand{\argmin}{\mathop{\mathrm{arg\,min}}}

\newcolumntype{G}{>{\centering\columncolor{gray!20!white}}p{0.2\textwidth}}
\newcolumntype{C}{>{\centering\arraybackslash}p{0.15\textwidth}}
\newcolumntype{c}{>{\centering\arraybackslash}p{0.10\textwidth}}

\def\d{\delta}

\def\xinid{\xi_{n,i}^\delta}
\def\zetanid{\zeta_{n,i}^\delta}
\def\znid{z_{n,i}^\delta}
\def\munid{\mu_{n,i}^\delta}
\def\snid{\sigma_{n,i}^\delta}

\def\lnid{\lambda_{n,i}^\delta}

\def\D{\mathcal D}

\def\X{\mathcal X}

\usepackage{lineno,hyperref}

\title{Acousto-electric tomography by the convergence of Kaczamrz two-point gradient-$\Theta$ method}
\author{ Kai Zhu\quad Jijun Liu\quad Min Zhong \footnote{Corresponding author:  Min Zhong, email: min.zhong@seu.edu.cn}  \\
School of Mathematics, Southeast University, Nanjing  210096, P. R. China\\
Nanjing Center for Applied Mathematics, Nanjing 211135, P. R. China
}

\begin{document}

\begin{abstract}
We study the numerical reconstruction problem in acousto-electric tomography
(AET) of recovering the conductivity distribution in a bounded domain from multiple interior power density data. The Two-Point-Gradient-$\Theta$ (TPG-$\Theta$) in Kaczmarz type is proposed, with a general convex penalty term $\Theta$, the algorithm can be utilized in AET problem for recovering sparse and discontinuous conductivity distributions. We establish the convergence of such iterative regularized method. Extensive numerical experiments are
presented to illustrate the feasibility and effectiveness of the proposed approach.\\

\end{abstract}
\maketitle

\section{Introduction}
The conductivity value varies widely with soft tissue types \cite{kr1989dielectrical,morimoto1993study} and its accurate imaging can provide valuable information about the physiological and pathological conditions of the tissue. This practical demand has driven the rapid development of medical imaging science research. The electrical impedance tomography (EIT) is an emerging technology that aims at reconstructing the conductivity distribution in a body from electrostatic measurements of voltages and currents on the surface of the body. The EIT is well known for being severely ill-posed problem \cite{mandache2001exponential}, it is necessary and challenging to propose a convergent and stable algorithm to reconstruct the conductivity. A novel idea of coupling EIT with different physical phenomenon has been promoted in the last decades. For example, EIT modulated by magnetic ultrasound waves leads to magnetic resonance EIT \cite{seo2011magnetic}, the EIT modulated ultrasound waves leads to acousto-electric tomography (AET), or equivalently, IAT \cite{gebauer2008impedance}. All the modalities give rise to additional interior information, and the availability of such interior information lead to a significant improvement of the conductivity reconstructions.

In this paper, we consider the AET, which is one of the most popular hybrid tomography imaging that has received increasing interests in the last decades \cite{ammari2008introduction,bal2013cauchy,widlak2012hybrid,zhang2004acousto}. Mathematically, AET consists in the determination of the spatially varying conductivity $\sigma>0$ in a bounded domain  $\Omega \subset \mathbb{R}^d$  from $I$ measurements of power densities $H_i(\sigma) = \sigma|\nabla u_i(\sigma)|^2 (i=0,1,\cdots,I-1)$  inside the $\Omega$ , resulting from $I$ different injected currents $f_i$. That is, the potential $u_i$ satisfied the elliptic equation.
    \begin{align}
    \label{eq:EIT}
                \left\{\begin{aligned}
                -\nabla \cdot\left(\sigma \nabla u_i\right) & =0,&  & \text { in } \Omega, \\
                        \sigma \frac{\partial u_i}{\partial \nu}  &=f_i, &  &\text { on } \partial \Omega ,
                \end{aligned}\right.
    \end{align}
where $\nu$ denotes the unit outward normal vector on $\partial \Omega$. 

There are some theoretical results available on the problem of estimating conductivity from power density functions (or some additional boundary information). The issues of uniqueness and stability have been studied intensively, see \cite{alberti2018lectures, bal2013hybrid} and references therein. Various numerical reconstructions in AET problem has been considered. In \cite{ammari2008electrical}, an algorithm was proposed for recovering conductivity from multiple power density measurements. Two optimal control formulations for reconstructing conductivity were introduced \cite{capdeboscq2009imaging}. The Levenberg-Marquardt iteration method \cite{2013TheLM} has been proposed, but only considered in a smooth Hilbert space $H^{\tau}(\Omega)$ with $\tau>d/2$ and data is noise free.  The explicit formulation of the reconstruction problems as regularized output least-squares problems was achieved \cite{adesokan2018fully}. A linearized reconstruction technique was developed in \cite{hoffmann2014iterative,yazdanian2021numerical} provided a general formulation of impedance tomography by the Picard and Newton iterative scheme. Recently, \cite{jin2022imaging} analyzed the problem in a view of learning, a deep neural network was considered which enjoyed remarkable robustness with respect to the presence of a large amount of data noise.

In recent years, $L^1$-liked or total variation(TV)-liked penalties are extremely popular for image processing, especially when the sought solution has special features such as sparsity or discontinuity. There are several works on nonlinear parameter identification problems in PDE with  $L^1$-liked or TV-liked penalty terms, e.g.,  EIT\cite{MR3767699}, quantitative photo-acoustic tomography\cite{MR3537006}, diffuse optical tomography\cite{MR4399528}, MAT-MI \cite{MR4593933}, acousto-electric tomography\cite{MR3916468}. In some of these works, iterative regularization methods with general convex penalty terms were utilized. In this paper, we propose a new Kaczmarz method for obtaining the stable approximation of AET conductivity. The nonlinear ill-posed operator equations are of the form
\begin{align*}
    H_i(\sigma) = y_i, \quad i=0,1,...,I-1 \,,
\end{align*}
in which $y_i = \sigma |\nabla u_i(\sigma)|^2$ and $u_i(\sigma)$ be the weak solution of \eqref{eq:EIT}. The operators $H_i:\mathcal{X}\rightarrow\mathcal{Y}$ are parameter-to-observation mappings between Hilbert space $\mathcal{X}$ and Banach space $\mathcal{Y}$. Instead of exact $y_i$, the noisy $y_i^\delta\in\mathcal{Y}$ are available, satisfying
\begin{align*}
\|y_i-y_i^\delta\|\leq\delta_i\,,\quad i=0,1,\cdots, I-1\,.
\end{align*}
A classical Landweber-Kaczmarz iteration with general convex penalty terms for was introduced in \cite{jin2013landweber}. Let $\Theta:\mathcal{X}\rightarrow(-\infty,\infty]$ be a proper, lower semi-continuous and uniformly convex functional, the method can be formulated as
\begin{align*}
&\xi_{n,i+1}^\delta = \xi_{n,i}^\delta-\mu_{n,i}^\delta H_i^\prime(\sigma_{n,i}^\delta)^*J_r^{\mathcal{Y}}\left(H_i(\sigma_{n,i}^\delta)-y_i^\delta\right)\,,\\
&\sigma_{n,i+1}^\delta = \textrm{arg}\min_{\sigma\in\mathcal{X}}\{\Theta(\sigma)-\langle\xi_{n,i+1}^\delta,\sigma\rangle\}\,,\quad i=0,1,\cdots, I-1\,.
\end{align*}
where $H_i'(\sigma)$ and $H_i'(\sigma)^*$ denote the Fr\'{e}chet derivative of $H_i$ at $\sigma$ and its adjoint. The $\mu_{n,i}^\delta$ is step size, $J_r^{\mathcal{Y}}:\mathcal{Y}\rightarrow \mathcal{Y}^*$ with $1<r<\infty$ denotes the duality mapping with gauge function $t\rightarrow t^{r-1}$. The advantage of this method is the freedom on the choice $\Theta$ so that it can be utilized in detecting the different features of the sought solution. However, the Landweber-typed methods are a type of slowly methods, by incorporating an extrapolation step into the iteration, we consider the two point gradient-$\Theta$ (TPG-$\Theta$) method as an acceleration 
\begin{equation*}
\begin{aligned}
&\zeta_{n,i}^\delta = \xi_{n,i}^\delta+\lambda_{n,i}^\delta(\xi_{n,i}^\delta-\xi_{n,i-1}^\delta)\,,\\
&z_{n,i}^\delta = \textrm{arg}\min_{\sigma\in\mathcal{X}}\{\Theta(\sigma)-\langle\zeta_{n,i}^\delta,\sigma\rangle\}\,,\\
&\xi_{n,i+1}^\delta = \zeta_{n,i}^\delta-\mu_{n,i}^\delta H_i'(z_{n,i}^\delta)^*J_r^{\mathcal{Y}_i}\left(H_i(z_{n,i}^\delta)-y_i^\delta\right)\,,\quad i=0,1,\cdots, I-1\,.
\end{aligned}
\end{equation*}
in which $\lambda_{n,i}^\delta$ is combination parameter. After an appropriate stopping criteria, the final iterate 
\begin{align*}
\sigma_{n,I}^\delta  = \textrm{arg}\min_{\sigma\in\mathcal{X}}\{\Theta(\sigma)-\langle\xi_{n,I}^\delta,\sigma\rangle\}
\end{align*}
will be regarded as the approximated conductivity. The TPG related methods have been confirmed its excellent acceleration in \cite{Ramlau2017,zhong2019regularization}, which discussed the single nonlinear operator in Hilbert space and Banach space respectively. We will provide detailed convergence analysis for such TPG-Kaczmarz type method and provide the numerical simulations to show its effectiveness and reliability.

The paper is organized as follows. In section 2, we give some properties of AET forward operator and convex analysis. In section 3, we formulate TPG-$\Theta$ method of Kaczmarz type and present the detailed convergence analysis. Finally in section 4, numerical simulations for AET problem are provided to test the performance of the method, geometrical and brain phantom are included in details, respectively.
\section{Preliminaries}
In this section, we introduce some necessary concepts and properties related to AET forward problem and convex analysis.

Let $\Omega$ be a nonempty, bounded, open, and connected set in $\mathbb{R}^d$, the boundary $\Gamma$ is Lipschitz continuous. The space $H_{\diamond}^1(\Omega)\subset H^1(\Omega)$ consists of functions with zero mean on the boundary, that is
\begin{align*}
H_{\diamond}^1(\Omega): = \left\{u\in H^1(\Omega):\int_\Gamma u \,ds = 0\right\}\,.
\end{align*} 
It is well known from \cite{attouch2014variational} that if a function $\nu\in H_{\diamond}^1(\Omega)$, it satisfies the Poincar\'{e} type inequality that $\|\nu\|_{H^1(\Omega)}\leq C\|\nabla \nu \|_{L^2(\Omega)}$.

For fixed $\lambda\in(0,1)$, define the set
$\mathcal{S} = \{\sigma:\lambda\leq\sigma\leq\frac{1}{\lambda},\,\textrm{a.e.} \,\textrm{in}\ \Omega\},$ assume the exact solution $\sigma^\dag\in\mathcal{S}$, which particularly implies $\sigma^\dagger \in L^\infty(\Omega)$.
Denote $(\mathcal{S}, \tau)$ as 
the set $\mathcal{S}$ endowed with a topology $\tau$.

Given $f\in H_{\diamond}^{-\frac{1}{2}}(\Gamma)$,which is the dual space of $H^{\frac{1}{2}}(\Gamma)$,  the forward AET problem considers the Neumann problem \eqref{eq:EIT}. It is well known from standard theory for elliptic PDE \cite{evans2022partial} that the Neumann problem has a unique weak solution $u(\sigma)\in V\equiv H_{\diamond}^1(\Omega)$ satisfying
\begin{align*}
(\sigma\nabla u,\nabla\psi) = \langle f, \psi\rangle\,,\quad \forall \psi\in V\,.
\end{align*}

Now we recall the regularity result for elliptic problems in $W^{1,q}(\Omega)$, see \cite{gallouet1999regularity} for detailed proof and \cite{1957Generalized,1963An} for related results.
\begin{theorem}\label{the}
        Let $\sigma \in \mathcal{S}$ and suppose $g \in L^{q}(\Omega), h \in L^{q}(\Omega)^{d}$
        and $f \in\left(W^{1-\frac{1}{q}, q}(\Gamma)\right)^{\prime}$ with $\int_{\Gamma} f \mathrm{~d} s+\int_{\Omega} g \mathrm{~d} x=0$.
        Then there exists a constant $Q>2$, such that for any $q \in(2, Q)$, the problem
        $$
        \left\{\begin{aligned}
        -\nabla \cdot(\sigma \nabla u) &=g+\nabla \cdot h, & & \text { in } \Omega, \\
        \sigma \frac{\partial u}{\partial \nu} &=f, & & \text { on } \Gamma,
        \end{aligned}\right.
        $$
        has a unique weak solution $u(\sigma) \in W^{1, q}(\Omega)$ satisfying
        $$
        \|u\|_{W^{1 , q}(\Omega)} \leqslant C\left(\|f\|_{\left(W^{1-\frac{1}{q} , q}(\Gamma)\right)^{\prime}}+\|g\|_{L^{q}(\Omega)}+\|h\|_{L^{q}(\Omega)^{d}}\right).
        $$
        The constant $Q=Q(\lambda, d)$ depends only on the domain $\Omega$, the spatial dimension $d$ and the constant $\lambda$, and the constant C depends only on $\Omega, \lambda, d$ and $q$.
\end{theorem}
In order to analyze the forward nonlinear map $\sigma\rightarrow H(\sigma)$, we first address the continuity of the solution operator $\sigma\rightarrow u(\sigma)$. This map is identical with that for EIT and has been extensively studied in various function spaces \cite{Chen1999An,2016THE,2012AnEIT}. Hence, we only list the results and the sketch of the proof, the details can be read in \cite{MR3916468}.
\begin{lemma}
    \label{lem:u_continuity}
Let $\{\sigma_k\}\subset\mathcal{S}$ with $\sigma_k\rightarrow\sigma^\dag$ in $L^1(\Omega)$, if $f\in \left(W^{1-\frac{1}{r},r}(\Gamma)\right)'$ for some $r>2$, then $u(\sigma_k)\rightarrow u(\sigma^\dag)$ in $W^{1,q}(\Omega)$ for any $q\in(2,\min(Q,r))$.
\end{lemma}
The next result gives the formula for the directional derivative $u'(\sigma)[\kappa]$ of $u(\sigma)$ at $\sigma$ on the direction $\kappa\in L^\infty(\Omega)$.
\begin{lemma}
    \label{lem:u_derivative}
For each $\sigma\in\mathcal{S}$, the direction derivatives $u'(\sigma)[\kappa]$ satisfies
\begin{align*}
(\sigma\nabla u'(\sigma)[\kappa],\nabla\varphi) = -(\kappa\nabla u(\sigma),\nabla\varphi)\,,\quad \forall \varphi\in V\,.
\end{align*}
If $f\in \left(W^{1-\frac{1}{r},r}(\Gamma)\right)'$ for some $r>2$, then for any $q\in(2,\min(Q,r))$ and $p'>\frac{q\min(Q,r)}{\min(Q,r)-q}$, then $u'(\sigma)[\cdot]: L^{p'}(\Omega)\rightarrow W^{1,q}(\Omega)$ is continuous, and $u(\sigma)$ is Fr\'{e}chet differentiable.
\end{lemma}

\begin{proof}
The expression of the directional derivative can be directly derived. If $f\in \left(W^{1-\frac{1}{r},r}(\Gamma)\right)'$, referring to Theorem \ref{the}, $u(\sigma)\in W^{1,q}(\Omega)$ for arbitrary $q\in(2,\min(Q,r))$. Thus we can choose any $q'\in(q,\min(Q,r))$ and $p'$ satisfying $(p')^{-1}+(q')^{-1} = q^{-1}$, such that 
\begin{align*}
\|\kappa\nabla u(\sigma)\|_{L^q(\Omega)}\leq \|\kappa\|_{L^{p'}(\Omega)}\|\nabla u(\sigma)\|_{L^{q'}(\Omega)}\,,
\end{align*}
and consequently, by Theorem \ref{the}
\begin{align}\label{eq:usigmakappa}
\|u'(\sigma)[\kappa]\|_{W^{1,q}(\Omega)}\leq c \|\kappa\|_{L^{p'}(\Omega)}\|\nabla u(\sigma)\|_{L^{q'}(\Omega)}\,.
\end{align}
This shows the boundness of $u'(\sigma)[\cdot]: L^{p'}(\Omega)\rightarrow W^{1,q}(\Omega)$, with arbitrary $p'>\frac{q\min(Q,r)}{\min(Q,r)-q}$.

Let $w(\sigma,\kappa) = u(\sigma+\kappa)-u(\sigma)-u'(\sigma)[\kappa]$, the residual satisfies
\begin{align*}
((\sigma+\kappa)\nabla w,\nabla\varphi) = -(\kappa \nabla u'(\sigma)[\kappa],\nabla\varphi),\quad \forall\varphi\in V\,.
\end{align*}
By the Theorem \ref{the} again, it follows that, for any $q\in(2,\min(Q,r))$,
\begin{align}\label{eq:w}
\|w(\sigma,\kappa)\|_{W^{1,q}(\Omega)}
\leq c\|\kappa\|_{L^{p'}(\Omega)}\|\nabla u'(\sigma)[\kappa]\|_{L^{q'}(\Omega)}\,
\leq c\|\kappa\|_{L^{p'}(\Omega)}^{1+p'/p''}\,,
\end{align}
in which $(p')^{-1}+(q')^{-1} = q^{-1}$, $p'$ is arbitrarily fixed by $p'>\frac{q\min(Q,r)}{\min(Q,r)-q}$ and $(p'')^{-1} + (q'')^{-1} = (q')^{-1}$. 
This shows the Fr\'{e}chet differentiability of $u:(\mathcal{S},L^{p'}(\Omega))\rightarrow W^{1,q}(\Omega)$.
\end{proof}

The following results give the continuity and differentiability of the forward map $H(\sigma):\sigma\rightarrow\sigma|\nabla u|^2$.

\begin{lemma}
        \label{lem:continuity}
        Let $\{ \sigma_k\} \subset \mathcal{S} $, $\sigma_k \rightarrow \sigma^\dagger$ in $L^1(\Omega)$, Then if $f \in (W^{1-\frac{1}{r},r}(\Gamma))^\prime$ for some $r>2$,
        Then $H(\sigma_k) \rightarrow H(\sigma^\dagger)$ in $L^{\frac{q}{2}}(\Omega)$ for arbitrary $q \in (2, \min(Q,r))$.
\end{lemma}

\begin{theorem}
The directional derivative $H'(\sigma)[\kappa]$ is given by
        \label{th:Hderivative}
        \begin{align*}
                H^\prime(\sigma)[\kappa] = \kappa | \nabla u(\sigma) |^2 + 2 \sigma \nabla u(\sigma) \cdot \nabla u^\prime(\sigma)[\kappa]\,.
        \end{align*}
If $f \in (W^{1-\frac{1}{r},r}(\Gamma))^\prime$ for some $r>2$, then for any $q\in(2,\min(Q,r))$ and $\kappa\in L^\infty(\Omega)$, $H'(\sigma)[\kappa]\in L^{q/2}(\Omega)$, and the adjoint $H'(\sigma)^*\omega$ is given by
\begin{align*}
H'(\sigma)^*\omega = |\nabla u|^2\omega +2\nabla u \cdot\nabla v\,,
\end{align*}
where $v = v(\omega)$ satisfies the
\begin{align}
    \label{eq:nu}
(\sigma\nabla v,\nabla\phi) = -(\sigma\omega\nabla u,\nabla\phi)\,,\quad\forall \phi\in V\,.
\end{align}
\end{theorem}

\begin{proof}
       The formal differentiation directly calculates by the chain rule 
        \begin{align*}
            H^\prime(\sigma)[\kappa] :=2 \sigma\nabla u^\prime(\sigma)[\kappa] \cdot \nabla u(\sigma)+\kappa |\nabla u(\sigma) |^2.
        \end{align*}

        Now we verify the boundness of $H^\prime(\sigma)[\kappa]$ in $L^{q/2}(\Omega)$. Referring to \eqref{eq:usigmakappa}, $u'(\sigma)[\kappa]\in W^{1,q}(\Omega)$ for arbitrary $q\in (2,\min(Q,r))$, thus
        \begin{align*}
        \|\kappa|\nabla u(\sigma)|^2\|_{L^{q/2}(\Omega)}
        &\leq \|\kappa\|_{L^{p'}(\Omega)}\|\nabla u(\sigma)\|^2_{L^{q'}(\Omega)}\,,
        \end{align*}
        in which the exponent $q'\in(q,\min(Q,r))$ and $p' = \frac{qq'}{2q'-q}$. In addition,
        \begin{align*}
        \|\sigma\nabla u'(\sigma)[\kappa]\cdot\nabla u(\sigma)\|_{L^{\frac{q}{2}}(\Omega)}&\leq c\|u'(\sigma)[\kappa]\|_{W^{1,q}(\Omega)}\|\nabla u(\sigma)\|_{L^q(\Omega)}\\
        &\leq c\|\kappa\|_{L^{p''}(\Omega)}\|\nabla u(\sigma)\|_{L^{q''}(\Omega)}\,,
        \end{align*}
        in which $(p'')^{-1}+(q'')^{-1} = q^{-1}$. Denote the Taylor Remainder
        \begin{align*}
            R(\sigma,\kappa) &= H(\sigma + \kappa) - H(\kappa) - H^\prime(\sigma)[\kappa]\\
            & = (\sigma+\kappa)|\nabla u(\sigma+\kappa)|^2-\sigma|\nabla u(\sigma)|^2-\kappa|\nabla u(\sigma)|^2-2\sigma\nabla u'(\sigma)[\kappa] \cdot \nabla u(\sigma)\,.
        \end{align*}
        Note that $ w(\sigma,k) = u(\sigma+\kappa)-u(\sigma)-u'(\sigma)[\kappa]$, yielding,
        \begin{align*}
        R(\sigma,\kappa) 
        & = (\sigma+\kappa)\left\vert\nabla(u(\sigma)+u'(\sigma)[\kappa]+w(\sigma,\kappa))\right\vert^2-\sigma|\nabla u(\sigma)|^2-\kappa|\nabla u(\sigma)|^2-2\sigma\nabla u'(\sigma)[\kappa] \cdot \nabla u(\sigma)\\
        & = (\sigma+\kappa)\nabla\left(w(\sigma,\kappa)+2u+2u'(\sigma)[\kappa]\right) \cdot \nabla w(\sigma,\kappa) + (\sigma+\kappa)|\nabla u'(\sigma)[\kappa]|^2+2\kappa\nabla u(\sigma)\cdot \nabla u'(\sigma)[\kappa]\,.
        \end{align*}
       We will estimate the $\|R(\sigma,\kappa)\|_{L^{q/2}(\Omega)}$ separately. First, referring to \eqref{eq:usigmakappa},
       \begin{align*}
       \|\kappa\nabla u(\sigma)\cdot \nabla u'(\sigma)[\kappa]\|_{L^{q/2}(\Omega)}\leq c\|\kappa\|_{L^{p'}(\Omega)}^2\|\nabla u(\sigma)\|_{L^{q'}(\Omega)}^2\,,
       \end{align*}
       in which the exponents $p',q'$ are defined as before. Then,
       \begin{align*}
       \|(\sigma+\kappa)|\nabla u'(\sigma)[\kappa]|^2\|_{L^{q/2}(\Omega)}&\leq
       \|\sigma+\kappa\|_{L^q(\Omega)}|\nabla u'(\sigma)[\kappa]|\|_{L^{q}(\Omega)}^2\\
       &\leq c\|\sigma+\kappa\|_{L^q(\Omega)}\|\kappa\|_{L^{p'}(\Omega)}^2\|\nabla u(\sigma)\|_{L^{q'}(\Omega)}^2\,.
       \end{align*}
       Finally, 
       \begin{align*}
       \|(\sigma+\kappa)\nabla\left(w(\sigma,\kappa)+2u+2u'(\sigma)[\kappa]\right)\cdot\nabla w(\sigma,\kappa)\|_{L^{q/2}(\Omega)}&\leq
       c\|\nabla(w+2u+2u'(\sigma)[\kappa]\|_{L^q(\Omega)}\|\nabla w\|_{L^q(\Omega)}\\
       &\leq c \|\kappa\|^{1+p'/p''}_{L^{p'}(\Omega)}\,,
       \end{align*}
       where the last inequality is from \eqref{eq:w}. This shows the Fr\'{e}chet differentiability of $H$.

       The adjoint $H'(\sigma)^*$ can be obtained by considering the weak formulations for $u'(\sigma)[\kappa]$ and $v$, see \cite{adesokan2018fully} for details.
\end{proof}

Next, we introduce some necessary concepts and properties related to Banach space and convex analysis, we refer to \cite{zalinescu2002convex} for more details.

 Let $\Theta:\mathcal{X}\rightarrow (-\infty, \infty]$ be continuous and convex functional, we call $\Theta$ is proper if its effective domain is nonempty. The conjugate $\Theta^*$ is defined by 
 \begin{align*}
 \Theta^{*}(\xi):=\sup _{x \in \mathcal{X}}\{\langle\xi, x\rangle-\Theta(x)\}, \quad \xi \in \mathcal{X}^*.
\end{align*}
 If $\Theta$ is proper, lower semi-continuous, and convex, the $\Theta^*$ shares the same properties. 
 
 Let $\partial \Theta(x)$ denote the subdifferential of $\Theta$ at $x \in \mathcal{X}$, let the Bregman distance induced by $\Theta$ at $x$ in the direction $\xi \in \partial \Theta(x)$ is defined by
 \begin{align*}
 D_{\xi} \Theta(\bar{x}, x):=\Theta(\bar{x})-\Theta(x)-\langle\xi, \bar{x}-x\rangle, \quad \forall \bar{x} \in \mathcal{X} .
\end{align*}
 Clearly, $D_{\xi} \Theta(\bar{x}, x) \geq 0$, and the following identity
 \begin{align}
 \label{eq:bregman distance}
 D_{\xi_{2}} \Theta\left(x, x_{2}\right)-D_{\xi_{1}} \Theta\left(x, x_{1}\right)=
 D_{\xi_{2}} \Theta\left(x_{1}, x_{2}\right)+\left\langle\xi_{2}-\xi_{1}, x_{1}-x\right\rangle
 \end{align}
 is valid for all $x \in \mathcal{D}(\Theta), x_{1}, x_{2} \in \mathcal{D}(\partial \Theta)$, and $\xi_{1} \in \partial \Theta\left(x_{1}\right), \xi_{2} \in \partial \Theta\left(x_{2}\right)$.
 
Assume the proper convex function $\Theta: \mathcal{X} \rightarrow(-\infty, \infty]$ is strongly convex,  then for arbitrary $x, \bar{x} \in \mathcal{D}(\partial \Theta), \xi \in \partial \Theta(x), \bar{\xi} \in$ $\partial \Theta(\bar{x})$, 
 it can be proved in \cite{zhong2019regularization} that
 \begin{align}\label{eq:leq}
 D_{\xi} \Theta(\bar{x}, x) \leq \frac{1}{4 c_{0}}\|\xi-\bar{\xi}\|^{2}.
 \end{align}
 Moreover, it follows from \cite{zalinescu2002convex}(Corollary 3.5.11) that $\mathcal{D}\left(\Theta^{*}\right)=\mathcal{X}^*, \Theta^{*}$ is Fr\'{e}chet differentiable and its gradient $\nabla \Theta^{*}: \mathcal{X}^* \rightarrow \mathcal{X}$ satisfies
 \begin{align*}
 \left\|\nabla \Theta^{*}\left(\xi_{1}\right)-\nabla \Theta^{*}\left(\xi_{2}\right)\right\| \leq \frac{\left\|\xi_{1}-\xi_{2}\right\|}{2 c_{0}}, \quad \forall \xi_{1}, \xi_{2} \in \mathcal{X}^*.
\end{align*}
 Consequently, it follows from \cite{zalinescu2002convex} that
 \begin{align*}
 x \in \nabla \Theta^{*}(\xi) \Longleftrightarrow \xi \in \partial \Theta(x) \Longleftrightarrow x=\arg \min _{z \in \mathcal{X}}\{\Theta(z)-\langle\xi, z\rangle\} \,.
\end{align*}

\section{Iterative TPG-$\Theta$ method in Kaczmarz type}
Instead of considering TPG-$\Theta$ method for solving the nonlinear equation $H(\sigma) = y$ directly, we consider a more general setup in which the single equation is extended to a power density system
\begin{align}\label{eq:eqs}
H_i(\sigma) = y_i\,,\quad i=0,1,\cdots I-1\,,
\end{align}
consisiting of $I$ equations, in which $y_i = \sigma|\nabla u_i|^2$, and $u_i(\sigma)$ be the weak solution of the boundary value problem \eqref{eq:EIT} with boundary value $f_i \in (W^{1-\frac{1}{r},r}(\partial \Omega))^\prime$. For each $i=0,1,...,I-1$, the power density operator $H_i:\mathcal{S}\subset\mathcal{X}\rightarrow\mathcal{Y}_i$, parameter-to-observing mappings between a real Hilbert space $\X$ to a Banach space $\mathcal{Y}_i$. Such system arises naturally in many practical applications including AET with multiple exterior measurements.

By introducing
\begin{align*}
{\bf {H}}:(H_0,H_1,\cdots, H_{I-1}):\mathcal{X}^I\rightarrow \mathcal{Y}_0\times\mathcal{Y}_1\times\cdots\times\mathcal{Y}_{I-1}
\end{align*}
and ${\bf {y}} = (y_0,y_1,\cdots,y_{I-1})$. The system could be reformulated as a single equation ${\bf {H}}(\sigma) = {\bf {y}}$. However, it owns advantages to consider each equation separately and the memory consumption are significantly saved. We may work under the following conditions on $H_i$.
\begin{assumption}\label{ass1}
\begin{enumerate}
  \item There exists $\rho>0$, the ball $B_{2\rho}(\sigma_0)\subset\mathcal{S}$, and the equation  \eqref{eq:eqs} has a solution $\sigma^\dag\in \mathcal{D}(\Theta)$ such that
  \begin{align*}
  D_{\xi_0}\Theta(\sigma^\dag,\sigma_0)\leq c_0\rho^2\,.
  \end{align*}
  \item Each $H_i$ is Fr\'{e}chet differentiable on $B_{2\rho}(\sigma_0)$ and $\sigma\rightarrow H_i'(\sigma)$ is continuous in $B_{2\rho}(\sigma_0)$.
  \item For each $H_i$, there exists $0<\eta<1$ such that
  \begin{align*}
  &\|H_i'(\sigma)\|\leq C_H\,,\quad\sigma \in B_{3\rho}(\sigma_0)\,,\\
  &\|H_i(\sigma)-H_i(\bar \sigma)-H_i'(\bar \sigma)(\sigma-\bar \sigma)\|\leq \eta\|H_i(\sigma)-H_i(\bar \sigma)\|\,,\quad \sigma,\bar \sigma\in B_{3\rho}(\sigma_0)\,.
  \end{align*}
\end{enumerate}
\end{assumption}

In practical application, instead of $y_i$ we only have noisy data $y_i^\delta$ satisfying
\begin{align*}
\|y_i-y_i^\delta\|\leq\delta\,,\quad i=0,1,\cdots,I-1\,,
\end{align*}
we will use $y_i^\delta$, $i=0,1,\cdots,I-1$ to reconstruct the solution of equation \eqref{eq:eqs}. In the following manuscript, We will assume $\mathcal{X} = L^2(\Omega)$, each $\mathcal{Y}_i$ be the Banach space $L^{q/2}(\Omega)$ with $q\in(2,\min(Q,r))$. Such spaces are uniformly smooth, thus the duality mappings $J_{q/2}:\mathcal{Y}_i\rightarrow \mathcal{Y}_i^*$ are single valued and continuous. The TPG-$\Theta$ method in Kaczmarz type is proposed in the following Algorithm \ref{A1}.

\begin{algorithm}[h]
        \begin{algorithmic}
        \STATE   \textbf{i)} Pick the initial guess $\xi_0^\delta = \xi_0$ and $\sigma_0^\delta = \sigma_0 = \textrm{arg}\min_{\sigma\in L^2(\Omega)}\left\{\Theta(\sigma)-\langle\xi_0^\delta,\sigma\rangle\right\}\,.$
        \STATE  \textbf{ii)} Assume $\xi_n^\delta$ and $\sigma_n^\delta$ are defined for fixed $n\geq 0$. Set the initials $\xi_{n,-1}^\delta$ and $\xi_{n,0}^\delta$.
        \STATE  \qquad \textbf{For $i=0,1,\cdots I-1$,}
                \STATE  \qquad \textbf{choosing} the combination parameter $\lambda_{n,i}^\delta$, and define   
                                \begin{align*}
                                        \zetanid &= \xinid + \lnid(\xinid - \xi_{n,i-1}^\delta), \\
                                        \znid &= \argmin_{\sigma \in L^2(\Omega)} \{ \Theta(\sigma) - \left\langle \zetanid, \sigma \right\rangle \}\,,
                                \end{align*}
        \STATE \qquad \textbf{calculating} the residual  $r_{n,i}^\d=H_i(\znid) - y_i^\delta$,
        \STATE \qquad \textbf{calculating} the step size
                        \begin{align}\label{eq:mu}
                                \munid = \left\{
                                        \begin{array}{ll}
                                        \min\left\{\frac{\bar\mu_0\|r_{n,i}^\d\|^{2(q/2-1)}}{\|H^\prime_i(\znid)^*J_{q/2}(r_{n,i}^\d)\|^2},\bar\mu_1\right\}\| r_{n,i}^\delta \Vert^{2-q/2}, & \hbox{if $\|r_{n,i}^\d\|>\tau\d$;} \\
                                        0, & \hbox{if $\|r_{n,i}^\d\|\leq\tau\d$,}
                                        \end{array}
                                        \right.
                        \end{align}
        \STATE \qquad\qquad               for some constant $\bar\mu_0$ and $\bar\mu_1$, and $\tau>1$.
        \STATE \qquad \textbf{Define} $\xi_{n,i+1}^\delta = \zeta_{n,i}^\delta - \mu_{n,i}^\delta H^\prime_i(\znid)^*J_{q/2}(r_{n,i}^\delta)$.
        \STATE \qquad \textbf{Calculate} $\sigma_{n,i+1}^\d = \argmin_{\sigma \in \mathcal{X}} \{ \Theta(\sigma) - \left\langle \xi_{n,i+1}^\delta , \sigma \right\rangle \}$.
        \STATE \qquad \textbf{End}
        \STATE \qquad \textbf{Define} $\xi_{n+1}^\delta = \xi_{n,I}^\delta$ and $\sigma_{n+1}^\delta = \sigma_{n,I}^\d$,
        \STATE   \textbf{iii)} Let $n_\delta$ be the first integer s.t.,
        \begin{align}\label{eq:stop}
                \mu_{n_\delta,i}^\delta=0,\quad \text{for all} \quad i=0,1,...,I-1\,.
        \end{align}
        \STATE \textbf{iv)} Use $\sigma_{n_\delta}^\delta$ be the approximate solution.
        \end{algorithmic}
        \caption{TPG-$\Theta$ method in Kaczmarz type. \label{A1}}
\end{algorithm}
\subsection{Convergence}
To derive the convergence analysis, we are going to show that, for any solution $\hat\sigma$ of \eqref{eq:eqs} in $B_{2\rho}(\sigma_0)\bigcap\D(\Theta)$, the Bregman distance $D_{\xi_n^\d}\Theta(\hat\sigma,\sigma_n^\d)$ is monotonically decreasing with respect to $n$, provided that $0\leq n< n_\d$. To this end, introduce
\begin{align*}
\triangle_{n,i+1}^\d = D_{\xi_{n,i+1}^\d}\Theta(\hat\sigma, \sigma_{n,i+1}^\d)-D_{\xinid}\Theta(\hat \sigma,\snid)\,.
\end{align*}

We divide the $\triangle_{n,i+1}^\d$ into two parts
\begin{align*}
\triangle_{n,i+1}^\d &= D_{\xi_{n,i+1}^\d}\Theta(\hat\sigma, \sigma_{n,i+1}^\d) - D_{\zeta_{n,i}^\d}\Theta(\hat\sigma, z_{n,i}^\d) +  D_{\zeta_{n,i}^\d}\Theta(\hat\sigma, z_{n,i}^\d)  - D_{\xi_{n,i}^\d}\Theta(\hat\sigma, \sigma_{n,i}^\d)\\
&: = P_1 + P_2\,.
\end{align*}
Utilizing the equality \eqref{eq:bregman distance} and \eqref{eq:leq}, the first part can be estimated by
\begin{align*}
P_1 &= D_{\xi_{n,i+1}^\d}\Theta(z_{n,i}^\d,\sigma_{n,i+1}^\d) + \langle\xi_{n,i+1}^\d-\zeta_{n,i}^\d,\znid-\hat\sigma\rangle\\
&\leq \frac{1}{4c_0}\|\xi_{n,i+1}^\d - \zetanid\|^2 + \langle\xi_{n,i+1}^\d-\zeta_{n,i}^\d,\znid-\hat\sigma\rangle\,.
\end{align*}
Recalling the definition of $\xi_{n,i+1}^\d$, and also by the definition of $\munid$, one can see
\begin{align*}
\frac{1}{4c_0}\|\xi_{n,i+1}^\d - \zetanid\|^2 \leq \frac{\bar\mu_0}{4c_0}\munid\|r_{n,i}^\d\|^{q/2}\,.
\end{align*}
In addition, if $\znid\in B_{3\rho}(\sigma_0)$, the Assumption \ref{ass1}(iii) can be utilized to derive
\begin{align*}
\langle\xi_{n,i+1}^\d-\zeta_{n,i}^\d,\znid-\hat\sigma\rangle& = -\munid\langle J_{q/2}(r_{n,i}^\d), y_i^\d-H_i(\znid)-H_i'(\znid)(\hat \sigma-\znid)\rangle - \munid\|r_{n,i}^\d\|^{q/2}\\
&\leq\munid\|r_{n,i}^\d\|^{q/2-1}\left(\d + \| H_i(\hat \sigma)-H_i(\znid)-H_i'(\znid)(\hat \sigma-\znid)\|\right)-\munid\|r_{n,i}^\d\|^{q/2}\\
&\leq (1+\eta)\munid\|r_{n,i}^\d\|^{q/2-1}\d -(1-\eta)\munid\|r_{n,i}^\d\|^{q/2}\,,
\end{align*}
and since $\delta < \|r_{n,i}^\d\|/\tau$, it follows that,
\begin{align*}
P_1< -(1-\eta-\frac{1+\eta}{\tau}-\frac{\bar\mu_0}{4c_0})\munid\|r_{n,i}^\d\|^{q/2}\,.
\end{align*}
For the second term, applying \eqref{eq:bregman distance} and \eqref{eq:leq} again,
\begin{align*}
P_2 &= D_{\zetanid}\Theta(\sigma_{n,i}^\d,\znid) + \langle\zetanid-\xinid,\snid-\hat\sigma\rangle\\
&\leq \frac{1}{4c_0}(\lnid)^2\|\xinid-\xi_{n,i-1}^\d\|^2 + \langle\zetanid-\xinid,\snid-\hat\sigma\rangle\,.
\end{align*}
Since,
\begin{align*}
\langle\zetanid-\xinid,\snid-\hat\sigma\rangle&= \lnid\langle\xinid-\xi_{n,i-1}^\d,\snid-\hat\sigma\rangle\\
& = \lnid\left(D_{\xinid}\Theta(\hat\sigma,\snid)-D_{\xi_{n,i-1}^\d}\Theta(\hat\sigma,\sigma_{n,i-1}^\d)
+D_{\xi_{n,i-1}^\d}\Theta(\snid,\sigma_{n,i-1}^\d)\right)\\
& \leq \lnid\triangle_{n,i}+\frac{\lnid}{4c_0}\|\xinid-\xi_{n,i-1}^\d\|^2\,,
\end{align*}
and consequently
\begin{align}\label{eq:P2}
P_2\leq \lnid\triangle_{n,i}+\frac{(\lnid)^2+\lnid}{4c_0}\|\xinid-\xi_{n,i-1}^\d\|^2\,.
\end{align}
Combining the estimates of the two parts, we have the following lemma.
\begin{lemma}\label{lem1}
Let the Assumption \ref{ass1} is satisfied. Assume the $\tau>1$ and $\bar\mu_0>0$ are suitably chosen such that
\begin{align*}
c_1: = 1-\eta-\frac{(1+\eta)}{\tau}-\frac{\bar \mu_0}{4c_0}>0\,.
\end{align*}
Then, for any solution $\hat\sigma\in B_{2\rho}(\sigma_0)\bigcap\D(\Theta)$, for $n\geq 0$ and $i=0,1,\cdots, I-1$, if $\znid\in B_{3\rho}(\sigma_0)$, there holds,
\begin{align*}
\triangle_{n,i+1}< \lnid\triangle_{n,i}+\frac{(\lnid)^2+\lnid}{4c_0}\|\xinid-\xi_{n,i-1}^\d\|^2-c_1\munid\|r_{n,i}^\d\|^{q/2}\,.
\end{align*}
\end{lemma}

Now it is necessary to discuss the condition for positive combination parameter $\lnid$. Assume it satisfies
\begin{align}\label{eq:choiceadd}
\frac{(\lnid)^2+\lnid}{4c_0}\|\xinid-\xi_{n,i-1}^\d\|^2\leq c_0\rho^2\,,
\end{align}
and
\begin{align}\label{eq:choice}
\frac{(\lnid)^2+\lnid}{4c_0}\|\xinid-\xi_{n,i-1}^\d\|^2-\frac{c_1}{\nu}\munid\|r_{n,i}^\d\|^{q/2}\leq 0,
\end{align}
where $\nu>1$ is a constant independent of $\delta$ and $n$. In order to guarantee \eqref{eq:choice}, it is obvious that we only need consider the situation when $\|H_i(\znid)-y_i^\d\|>\tau\d$, since $\lnid$ will be forced to $0$ because $\munid = 0$ when $\|H_i(\znid)-y_i^\d\|\leq\tau\d$. Since $\munid \geq \min\left\{\frac{\bar\mu_0}{C_H^2},\bar\mu_1\right\}\|r_{n,i}^\d\|^{2-q/2}\,,$ it is sufficient to demand
\begin{align*}
\frac{(\lnid)^2+\lnid}{4c_0}\|\xinid-\xi_{n,i-1}^\d\|^2&\leq \frac{c_1}{\nu}\min\left\{\frac{\bar\mu_0}{C_H^2},\bar\mu_1\right\}\tau^2\d^2:= M\tau^2\d^2\,,
\end{align*}
which leads to the choice
\begin{align}\label{eq:choice3}
\lnid = \min\left\{-\frac{1}{2}+\sqrt{\frac{1}{4}+\frac{4c_0M\tau^2\d^2}{\|\xinid-\xi_{n,i-1}^\d\|^2}},\frac{n}{n+\alpha}\right\}\,,\quad\alpha\geq 3\,.
\end{align}
Note that in the above formula for $\lnid$, inside the “min”
the second argument is taken to be $n/(n+\alpha)$, which is the combination parameter utilized in Nesterov acceleration. 
\begin{proposition}\label{prop:1}
If the combination parameter $\{\lnid\}$ are chosen to satisfy \eqref{eq:choiceadd} and \eqref{eq:choice}, then for all $n\geq 0$ and $i=0,1,\cdots, I$, we have
\begin{align*}
\znid\in  B_{3\rho}(\sigma_0),\quad \snid\in B_{2\rho}(\sigma_0)\,.
\end{align*}
 for arbitrary solution $\hat \sigma\in B_{2\rho}(\sigma_0)\bigcap \D(\Theta)$, we have $\triangle_{n,i}\leq 0$, and particularly
\begin{align*}
D_{\xi_{n+1}^\d}\Theta(\hat\sigma, \sigma_{n+1}^\d)\leq D_{\xi_{n}^\d}\Theta(\hat \sigma, \sigma_{n}^\d)\,.
\end{align*}
Let $n_\d$ be chosen according to \eqref{eq:stop}, then $n_\d$ must be a finite number.
\end{proposition}

\begin{proof}
The conclusion is trivial for $n=0$ and $i=0$. Suppose the conclusion is hold for $\{\znid\}$, $\{\snid\}$ and $\triangle_{n,i}$ for arbitrary $0\leq n\leq m$ and $i=0,1,\cdots, I$. Then, when $n=m+1$, since $\xi_{m+1,-1}^\d=\xi_{m,I-1}^\d, \xi_{m+1,0}^\delta = \xi_{m,I}^\d$ and $\sigma_{m+1,-1}^\d  = \sigma_{m,I-1}^\d \in  B_{2\rho}(\sigma_0) \, , \sigma_{m+1,0}^\delta = \sigma_{m,I}^\d \in B_{2\rho}(\sigma_0)$, this leads $\triangle_{m+1,0}= \triangle_{m,I} \leq 0$. Combining \eqref{eq:choice}, the estimate
\begin{align*}
\triangle_{m+1,i+1}< \lambda_{m+1,i}^\d\triangle_{m+1,i}-(1-\frac{1}{\nu})c_1\mu_{m+1,i}^\d\|r_{m+1,i}^\d\|^{q/2}
\end{align*}
provides $\triangle_{m+1,i+1}\leq 0$. Therefore, taking $\hat\sigma = \sigma^\dag$,
\begin{align*}
D_{\xi_{m+1,i+1}^\d}\Theta(\sigma^\dag, \sigma_{m+1,i+1}^\d)\leq D_{\xi_{m+1,0}^\d}\Theta(\sigma^\dag, \sigma_{m+1,0}^\d)\leq\cdots\leq D_{\xi_{0,0}^\d}\Theta(\sigma^\dag, \sigma_{0,0}^\d)\leq c_0\rho^2\,,
\end{align*}
this implies $\|\sigma^\dag-\sigma^\delta_{m+1,i+1}\|\leq\rho$, this combines $\sigma^\dag\in B_\rho(\sigma_0)$ implies $\sigma^\delta_{m+1,i+1}\in B_{2\rho}(\sigma_0)$.

Now we may recall the estimate for $P_2$, combining with \eqref{eq:choiceadd}, yielding
\begin{align*}
D_{\zeta_{m+1,i+1}^\d}\Theta(\sigma^\dag,z_{m+1,i+1}^\d)&\leq D_{\xi_{m+1,i+1}^\d}\Theta(\sigma^\dag,\sigma_{m+1,i+1}^\d) +\lambda_{m+1,i+1}^\d\triangle_{m+1,i+1}+c_0\rho^2\leq 2c_0\rho^2\,.
\end{align*}
This consequently gives $z_{m+1,i+1}^\d\in B_{3\rho}(x_0)$. The conclusions are valid.

Finally, we will show $n_\d$ is a finite number. Recalling that
\begin{align*}
\triangle_{n,i+1}\leq -(1-\frac{1}{\nu})c_1\mu_{n,i}^\d\|r_{n,i}^\d\|^{q/2}\,,\quad \forall \, 0\leq n< n_\d\,,
\end{align*}
this implies
\begin{align}\label{eq:moreover}
(1-\frac{1}{\nu})c_1\sum_{i=0}^{I-1}\mu_{n,i}^\d\|r_{n,i}^\d\|^{q/2}\leq D_{\xi_{n,0}^\delta}\Theta(\sigma^\dag,\sigma_{n,0}^\d)-D_{\xi_{n,I}^\d}\Theta(\sigma^\dag,\sigma_{n,I}^\d)\,.
\end{align}
According to the definition of $n_\d$, for arbitrary $n<n_\delta$, there is at least one index $i_n\in\{0,1,\cdots, I-1\}$, such that $\|H_{i_n}(z_{n,i_n}^\d)-y_{i_n}^\d\|>\tau\delta$. Consequently,
\begin{align*}
\sum_{i=0}^{I-1}\mu_{n,i}^\d\|r_{n,i}^\d\|^{q/2}\geq \min\left\{\frac{\bar\mu_0}{C_H^2},\bar\mu_1\right\}\|r_{n,i_n}^\delta\|^{2}>C\tau^2\delta^2\,.
\end{align*}
Summing \eqref{eq:moreover} for $n=0$ to $n=n_\d-1$ and using the above inequality, the $n_\delta$ is finite.
\end{proof}

In order to establish the regularization property of the method, we need to consider the noise-free case. By dropping the superscript $\delta$ in all the quantities involved in Algorithm \ref{A1},it leads to noise-free TPG-$\Theta$ method with Kaczmarz type. We have the following lemma.
\begin{lemma}\label{exact}
 Let the Assumption \ref{ass1} is satisfied. Assume the $\bar\mu_0>0$ is suitably chosen such that
\begin{align*}
c_2: = 1-\eta-\frac{\bar \mu_0}{4c_0}>0\,,
\end{align*}
assume the combination parameter $\lambda_{n,i}$ is chosen satisfying  \eqref{eq:choiceadd} and \eqref{eq:choice} with $c_2$ replacing $c_1$. Then, for all $n\geq 0$ and $i=0,1,\cdots, I$, we have  $z_{n,i}\in B_{3\rho}(\sigma_0)$, $\sigma_{n,i}\in B_{2\rho}(\sigma_0)$, and for any solution $\hat\sigma\in B_{2\rho}(\sigma_0)\bigcap\D(\Theta)$, we have
\begin{align}
&D_{\xi_{n+1}}\Theta(\hat\sigma,\sigma_{n+1})\leq D_{\xi_n}\Theta(\hat\sigma,\sigma_n)\,,\label{eq:decrease2}\\
& (1-\frac{1}{\nu})c_2\sum_{i=0}^{I-1}\mu_{n,i}\|r_{n,i}\|^{\frac{q}{2}}\leq D_{\xi_{n,0}}\Theta(\hat\sigma,\sigma_{n,0})-D_{\xi_{n,I}}\Theta(\hat\sigma,\sigma_{n,I})\,.\label{eq:add2}
\end{align}
In addition, $\lim_{n\rightarrow\infty}\sum_{i=0}^{I-1}\|r_{n,i}\|^2=0\,.$

\end{lemma}

Now we are ready to establish convergence result for the two point gradient-$\Theta$ method in Kaczmarz type. The following proposition is important.
\begin{proposition}\label{useful}
Let Assumption \ref{ass1} is satisfied, let $\{\xi_n,\sigma_n\}$ satisfy
\begin{enumerate}
    \item $\sigma_n\in B_{2\rho}(\sigma_0)$, and $\xi_n\in\partial\Theta(\sigma_0)$ for all $n\geq 0$,
    \item $\{D_{\xi_n}\Theta(\hat\sigma,\sigma_n)\}$ is monotonically decreasing for arbitrary solution $\hat\sigma\in B_{2\rho}(\sigma_0)\bigcap\D(\Theta)$,
    \item $\lim_{n\rightarrow\infty}\|H_i(\sigma_n)-y_i\| = 0$,
    \item there exists a subsequence $\{n_k\}$ with $n_k\rightarrow\infty$ such that for arbitrary solution $\hat\sigma\in B_{2\rho}(\sigma_0)\bigcap\D(\Theta)$, there holds
    \begin{align*}
    \lim_{\ell\rightarrow\infty}\sup_{k\geq\ell}\vert\langle\xi_{n_k}-\xi_{n_\ell},\sigma_{n_k}-\hat\sigma\rangle\vert = 0.
    \end{align*}
\end{enumerate}
Then, there exists a solution $\sigma^*\in \X$ such that
\begin{align*}
\lim_{n\rightarrow\infty}\|\sigma_n-\sigma^*\| = 0\,,\quad \lim_{n\rightarrow\infty}D_{\xi_n}\Theta(\sigma^*,\sigma_n) = 0\,.
\end{align*}
\end{proposition}

The (i) and (ii) in Proposition \ref{useful} are already satisfied. Now we examine the (iii). Since $\|H_i(\sigma_n)-y_i\|$ can be splited into
\begin{align*}
\|H_i(\sigma_n)-y_i\|\leq \|H_i(z_{n,i})-y_i\| + \sum_{j=0}^i\|H_i(\sigma_{n,j})- H_i(z_{n,j})\| + \sum_{j=0}^{i-1}\|H_i(\sigma_{n,j+1})- H_i(z_{n,j})\|.
\end{align*}
Since by \eqref{eq:moreover}, and according to the choice strategy for $\lambda_{n,i}$ and $\mu_{n,i}$, 
\begin{align*}
\|z_{n,i}-\sigma_{n,i}\| &= \|\nabla\Theta^*(\zeta_{n,i})-\nabla\Theta^*(\xi_{n,i})\|\leq \frac{1}{2c_0}\|\zeta_{n,i}-\xi_{n,i}\|\\
&\leq \frac{\lambda_{n,i}}{2c_0}\|\xi_{n,i}-\xi_{n,i-1}\|\leq\frac{\sqrt{c_2\bar\mu_1}}{\sqrt{c_0\nu}}\|r_{n,i}\|\,.
\end{align*}

In addition, recalling the boundness of $H_i'$, we have
\begin{align*}
\|\sigma_{n,j+1}-z_{n,j}\|  &= \|\nabla\Theta^*(\xi_{n,j+1})-\nabla\Theta^*(\zeta_{n,j})\|\\
&\leq \frac{1}{2c_0}\mu_{n,j}\|H_j'(z_{n,j})^*J_{q/2}(r_{n,j})\|\leq \frac{C_H\bar\mu_1}{2c_0}\|r_{n,j}\| \,.
\end{align*}
Therefore, there exists the constant $C = 2\max\{\frac{C_H^2\bar\mu_1}{2(1-\eta)c_0},1,\frac{C_H\sqrt{c_2\bar\mu_1}}{(1-\eta)\sqrt{c_0\nu}}\}$, such that
\begin{align}\label{eq:rni}
\|H_i(\sigma_n)-y_i\|\leq C\sum_{j=0}^i\|r_{n,j}\| \rightarrow 0\,,\quad\textrm{as}\ n\rightarrow\infty\,.
\end{align}

Finally, in order to derive the convergence result, it is necessary to verify the (iv) in Proposition \ref{useful}. To this end, let
\begin{align}\label{eq:Rn}
R_n: = \sum_{i=0}^{I-1}\|r_{n,i}\|^2\,,
\end{align}
it follows from Lemma \ref{exact} that $\lim_{n\rightarrow\infty}R_n = 0$. Moreover, if $R_n = 0$ for some integer $n$, then $r_{n,i} = 0$ for all $i=0,1,\cdots, I-1$. The choice strategy for $\lambda_{n,i}$ provides that $\lambda_{n,i}(\xi_{n,i}-\xi_{n,i-1}) = 0$ and thus $\zeta_{n,i} = \xi_{n,i}$, then the choice strategy $\mu_{n,i} = 0$ yields $\xi_{n,i+1} = \zeta_{n,i} = \xi_{n,i}$ for all $i=0,1,\cdots, I-1$. This means the algorithm stop updating. It shows that
\begin{align}\label{eq:Rn2}
R_n = 0\ \textrm{for some}\ n\Rightarrow R_m = 0\ \textrm{for all}\ m\geq n\,. 
\end{align}
In view of \eqref{eq:Rn} and \eqref{eq:Rn2}, we can introduce a subsequence $\{n_k\}$ by letting $n_0 = 0$ and $n_k$ be the first integer satisfying 
\begin{align*}
n_k\geq n_{k-1}+1 \, \textrm{and}\ R_{n_k} \leq R_{n_{k-1}}\,.
\end{align*}
For such chosen strictly increasing sequence $\{n_k\}$, it is obvious 
\begin{align*}
R_{n_k}\leq R_n\,,\quad 0\leq n<n_k\,.
\end{align*}
We firstly consider that
\begin{align*}
\langle\xi_{n+1}-\xi_n,\sigma_{n_k}-\hat\sigma\rangle & = \sum_{i=0}^{I-1}\langle\xi_{n,i+1}-\xi_{n,i},\sigma_{n_k}-\hat\sigma\rangle\\
& = \sum_{i=0}^{I-1}\langle \lambda_{n,i}(\xi_{n,i}-\xi_{n,i-1})-\mu_{n,i}H_i'(z_{n,i})^*J_{q/2}(r_{n,i}),\sigma_{n_k}-\hat\sigma\rangle\,.
\end{align*}
Therefore, by using the property of duality mapping, we have
\begin{align*}
\vert\langle\xi_{n+1}-\xi_n,\sigma_{n_k}-\hat\sigma\rangle\vert&\leq \sum_{i=0}^{I-1}\lambda_{n,i}\|\xi_{n,i}-\xi_{n,i-1}\|\|\sigma_{n_k}-\hat\sigma\|\\
&\quad+ \sum_{i=0}^{I-1}\mu_{n,i}\|r_{n,i}\|^{q/2-1}\|H_i'(z_{n,i})(\sigma_{n_k}-\hat\sigma)\|\,.
\end{align*}
The first term can be estimated by
\begin{align}\label{eq:first}
\sum_{i=0}^{I-1}\lambda_{n,i}\|\xi_{n,i}-\xi_{n,i-1}\|\|\sigma_{n_k}-\hat\sigma\|\leq 4\rho\sum_{i=0}^{I-1}\lambda_{n,i}\|\xi_{n,i}-\xi_{n,i-1}\|\,.
\end{align}
By using the Assumption \ref{ass1} and Cauchy inequality, for $n<n_k$, the second term can be estimated by
\begin{align*}
 \sum_{i=0}^{I-1}\mu_{n,i}\|r_{n,i}\|^{q/2-1}\|H_i'(z_{n,i})(\sigma_{n_k}-\hat\sigma)\|
 &\leq \bar\mu_1\sum_{i=0}^{I-1}\|r_{n,i}\|(1+\eta)\left(\|H_i(z_{n,i})  - H_i(\sigma_{n_k})\| + \|H_i(z_{n,i})- y_i\|\right)\\
 &\leq (1+\eta)\bar\mu_1\sum_{i=0}^{I-1}\|r_{n,i}\|\|H_i(\sigma_{n_k})-y_i\|+2(1+\eta)\bar\mu_1R_n\\
 &\leq  (1+\eta)\bar\mu_1\sqrt{R_n}\left(\sum_{i=0}^{I-1}\|H_i(\sigma_{n_k})-y_i\|^2\right)^{1/2}+2(1+\eta)\bar\mu_1R_n\,.
\end{align*}
Recalling the \eqref{eq:rni}, it is obvious $\|H_i(\sigma_{n_k})-y_i\| \leq C\sqrt{I}\sqrt{R_{n_k}}\,,$ then,
\begin{align}\label{eq:second}
 \sum_{i=0}^{I-1}\mu_{n,i}\|r_{n,i}\|^{q/2-1}\|H_i'(z_{n,i})(\sigma_{n_k}-\hat\sigma)\|
 &\leq (1+\eta)\bar\mu_1CI\sqrt{R_n}\sqrt{R_{n_k}} + 2(1+\eta)\bar\mu_1R_n\nonumber\\
 &\leq \bar\mu_1(1+\eta)(CI+2)R_n\,.
\end{align}
The combination of \eqref{eq:first} and \eqref{eq:second} gives
\begin{align*}
\vert\langle\xi_{n_k}-\xi_{n_\ell},\sigma_{n_k}-\hat\sigma\vert&\leq \sum_{n=n_\ell}^{n_k-1}\vert\langle\xi_{n+1}-\xi_n,\sigma_{n_k}-\hat\sigma\vert\\
&\leq 4\rho \sum_{n=n_\ell}^{n_k-1}\sum_{i=0}^{I-1}\lambda_{n,i}\|\xi_{n,i}-\xi_{n,i-1}\|+ \bar\mu_1(1+\eta)(2+CI)\sum_{n=n_\ell}^{n_k-1}R_n\,.
\end{align*}
If we assume further that $\sum_{n=0}^\infty\sum_{i=0}^{I-1}\lambda_{n,i}\|\xi_{n,i}-\xi_{n,i-1}\|<\infty$, then the (iv) in Proposition \ref{useful} can be verified. We have the following convergence results for the TPG-$\Theta$ method in Kaczmarz type.
\begin{theorem}\label{th:convergence}
Let all the conditions in Lemma \ref{exact} hold, assume further 
\begin{align}\label{eq:lambda3}
\sum_{n=0}^\infty\sum_{i=0}^{I-1}\lambda_{n,i}\|\xi_{n,i}-\xi_{n,i-1}\|<\infty\,.
\end{align}
For the sequences $\{\xi_{n}\}$ and $\{\sigma_{n}\}$ defined by Algorithm \ref{A1} with exact data, there exists a solution $\sigma^*\in B_{2\rho}(x_0)\bigcap D(\Theta)$ such that
\begin{align*}
\lim_{n\rightarrow\infty}\|\sigma_n-\sigma^*\| = 0,\quad\lim_{n\rightarrow\infty} D_{\xi_n}\Theta(\sigma^*,\sigma_n) = 0\,.
\end{align*}
\end{theorem}

In order to prove the convergence for TPG-$\Theta$ method of Kaczmarz type in noisy case, we need the following stability results.
\begin{lemma}\label{lem:stability}
Let all the conditions in Lemma \ref{lem1} hold. Then for all $n\geq 0$, and $i=0,1,\cdots, I-1$, there hold
\begin{align*}
\xi_{n,i}^\d\rightarrow\xi_{n,i}, \ \zeta_{n,i}^\d\rightarrow\zeta_{n,i}, \ z_{n,i}^\d\rightarrow z_{n,i}, \ \sigma_{n,i}^\d\rightarrow \sigma_{n,i},\quad\textrm{as}\ \d\rightarrow 0\,.
\end{align*}
\end{lemma}
\begin{proof}
The results are trivial for $n=0$ and $i=0$. We assume that the results are true for some $n\geq 0$ and some index $i\in\{0,1,\cdots, I-1\}$ and show that they are valid for subscript $(n,i+1)$. To this end, we consider the two cases.

{\bf Case 1:} $\|r_{n,i}\| = 0$. In this case we have $\mu_{n,i} = 0$ and $\|r_{n,i}^\d\|\rightarrow 0$ due to the continuity of $H_i$ and the induction hypothesis $z_{n,i}^\d\rightarrow z_{n,i}$. Thus,
\begin{align*}
\xi_{n,i+1}^\d - \xi_{n,i+1} = \zeta_{n,i}^\d - \zeta_{n,i}-\mu_{n,i}^\d H_i'(z_{n,i}^\d)^* J_{q/2}(r_{n,i}^\d)\,,
\end{align*}
by the hypothesis for $\zeta_{n,i}^\d\rightarrow \zeta_{n,i}$, it is straightforward that
\begin{align*}
\|\xi_{n,i+1}^\d - \xi_{n,i+1}\|\leq \|\zeta_{n,i}^\d- \zeta_{n,i}\| + C_H\bar\mu_1\|r_{n,i}^\d\|\rightarrow 0 \,,\quad \textrm{as}\ \d\rightarrow 0\,.
\end{align*}
Consequently, the continuity of $\nabla\Theta^*$ yields $\sigma_{n,i+1}^\d \rightarrow \sigma_{n,i+1}$. In addition,
\begin{align*}
\zeta_{n,i+1}^\d = \xi_{n,i+1}^\d + \lambda_{n,i+1}^\d(\xi_{n,i+1}^\d - \xi_{n,i}^\d)\,,
\end{align*}
we may use $\lambda_{n,i+1}^\d\rightarrow \lambda_{n,i+1}$ to prove $\zeta_{n,i+1}^\d\rightarrow \zeta_{n,i+1}$ and immediately $z_{n,i+1}^\d \rightarrow z_{n,i+1}$ as $\d\rightarrow 0$.

{\bf Case 2:}  $\|r_{n,i}\| \neq 0$. In this case, we have $\|r_{n,i}^\d\|>\tau\d$ for small $\d$. The choice strategy for $\mu_{n,i}^\d$ provides the fact that $\mu_{n,i}^\d\rightarrow\mu_{n,i}$ as $\d\rightarrow 0$. By the continuity of $H_i$, $H_i'$, the duality mapping $J_{q/2}$ and $\nabla\Theta^*$, the conclusions are valid. 
\end{proof}

The final convergence result for TPG-$\Theta$ method in Kaczmarz type in noisy case can be proved similarly with Theorem 3.9 in \cite{jin2013landweber}.
\begin{theorem}\label{th:final}
Let all the conditions in Lemma \ref{lem1} hold, and the combination parameter $\lambda_{n,i}^\d$ are chosen satisfy \eqref{eq:choiceadd} and \eqref{eq:choice}. Then for $\{\xi_n^\d\}$ and $\{\sigma_n^\d\}$ defined by the TPG-$\Theta$ method of Kaczmarz type, there exists a solution $\sigma^*\in B_{2\rho}(x_0)\bigcap D(\Theta)$, such that
\begin{align*}
\lim_{\delta\rightarrow\infty} \sigma_{n_\delta}^\delta = \sigma^*,\quad \lim_{\delta\rightarrow\infty} D_{\xi_{n_\delta}^\delta}\Theta(\sigma^*,\sigma_{n_\delta}^\delta) = 0\,.
\end{align*}
\end{theorem}

\section{Numerical experiments}
The algorithm is implemented by Python using the DOLFIN 2019.2.0 (FEniCS) package
\cite{logg2010dolfin}, the standard solvers for PDE are formulated in the language of FEniCS with standard P1 finite elements. The meshes involved in generating simulated data and reconstruction obtained from the public software package \cite{geuzaine2009gmsh}. All meshes are generated by a circle shape, but with various refinement levels. In the experiments, the domain $\Omega$ is chosen by a disk in 2D in polar coordinates,
\begin{align*}
        \Omega = \{(r,\theta)\in [0,\frac{1}{2}] \times [0,2\pi)\}\,.
\end{align*}
 The $\mathcal{X} = L^2(\Omega)$ and $\mathcal{Y} = L^{q/2}(\Omega)$ with $q/2>1$ but close to $1$, which can be regarded as an approximation for $L^1$ fitting term. Its duality mapping $J_{q/2}: L^{q/2}\rightarrow L^{q/(q-2)}$ is defined by
 \begin{align*}
 J_{q/2}(\varphi) = |\varphi|^{q/2-1}\textrm{sign}(\varphi)\,,\quad\varphi\in L^{q/2}(\Omega)\,.
 \end{align*}
 We consider the boundary currents 
\begin{align*}
   \mathcal{F} = (f_1,f_2,f_3,f_4)=(x_1,x_2,\frac{x_1+x_2}{\sqrt{2}},\frac{x_1-x_2}{\sqrt{2}})\,.
\end{align*}
The simulated power density $\mathcal{H} = (H_1,H_2,H_3,H_4)$ corresponding to the boundary fluxes $\mathcal{F}$ are generated in a finer mesh $M_1$ (20054 nodes and 39748 triangles with mesh size $h=0.01$), and the reconstructions are performed on the coarser mesh $M_2$ (8272 nodes and 16140 triangles with mesh size $h=1/64$). The fine mesh can accurately resolve the forward problem and the coarse mesh can mitigate the so called inverse crime.

We consider two different phantoms: the geometrical shapes phantom and the brain phantom. The noisy data ${\bf H^\delta}$ is generated component-wise by
\begin{align}\label{eq:noise}
{\bf H^\delta} = {\bf H} + \delta_e\frac{\| {\bf H\Vert }}{\| {\bf e}\Vert }{\bf e},
\end{align}
where ${\bf e}$ is a vector whose entries satisfy a standard Gaussian distribution, $\delta_e$ is relative noise level, $\|{\bf H}\|$ denotes discrete $L^{q/2}$ norm in the finite element sense. It should be noted that the physically realistic absolute noise level is highly dependent on problem itself, this is due to the invert of the boundary data.  In our numerical experiments, the absolute noise level is $ \delta_e\|\bf{H}\|$. 

A key ingredient for numerical implementation is the determination of $\sigma = \nabla\Theta^*(\xi)$ for any given $\xi\in L^{2}(\Omega)$, which is equivalent to the minimization problem
\begin{align}\label{eq:calculate Fenchel}
\sigma = \textrm{arg}\min_{z\in L^2(\Omega)}\left\{\Theta(z)-\langle\xi,z\rangle\right\}\,.
\end{align}
 If the solution is sparse, we may choose
\begin{align*}
\Theta_{L^1}(\sigma) = \frac{1}{2\beta}\|\sigma\|_{L^2(\Omega)}^2 + \|\sigma\|_{L^1(\Omega)}\,,  \quad \beta>0\,,
\end{align*}
then the minimization problem \eqref{eq:calculate Fenchel} has the explicit form
\begin{align*}
\sigma = \beta\textrm{sign}(\xi)\max\{|\xi|-1,0\}\,.
\end{align*}
If the solution is piecewise constant, we may choose
\begin{align*}
\Theta_{TV}(\sigma) = \frac{1}{2\beta}\|\sigma\|_{L^2(\Omega)}^2 + \int_\Omega |D\sigma|\,,  \quad \beta>0\,,
\end{align*}
then the minimization problem \eqref{eq:calculate Fenchel} is equivalent to
\begin{align}\label{eq:ROF}
\sigma = \textrm{arg}\min_{z\in L^2(\Omega)}\left\{\frac{1}{2\beta}\|z-\beta\xi\|_{L^2(\Omega)}^2 +\int_\Omega |D z|\right\}\,,
\end{align}
which is a classical variational denoising problem. The fast iterative shrinkage-thresholding algorithm (FISTA) in \cite{beck2009fast,beck2009fastIEEE} is able to solve the \eqref{eq:ROF} in standard meshes, but can not be directly ported into trigonometric elements. Therefore, in our experiments, we use the primal-dual Newton method to solve such minimization problem\cite{MR1928831}.

\subsection{Numerical results for full data}
\begin{figure}[H]
\centering
 \begin{subfigure}{0.32\textwidth}
 \centering
\includegraphics[width = 5cm]{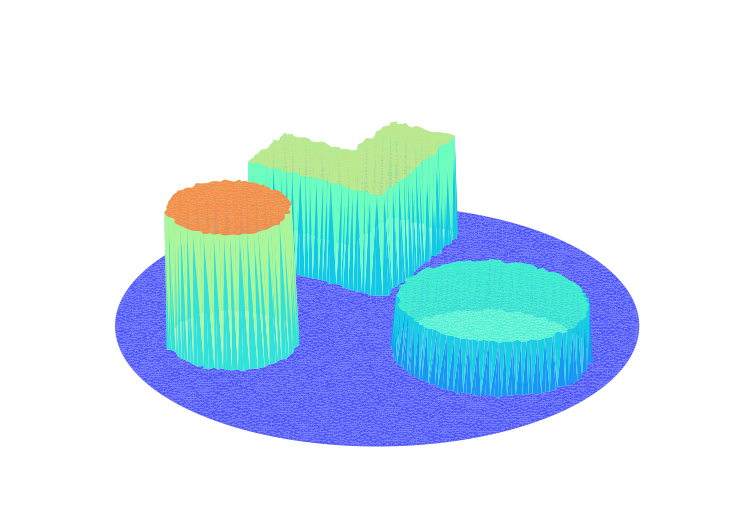}
\caption*{(a)}
\end{subfigure}
\begin{subfigure}{0.05\textwidth}
		\includegraphics[width=1.1cm,height = 3.8cm]{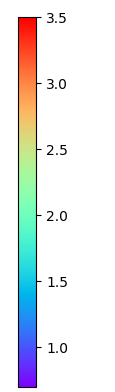}
		\caption*{ }
\end{subfigure}
 \begin{subfigure}{0.32\textwidth}
\includegraphics[width=3.6cm]{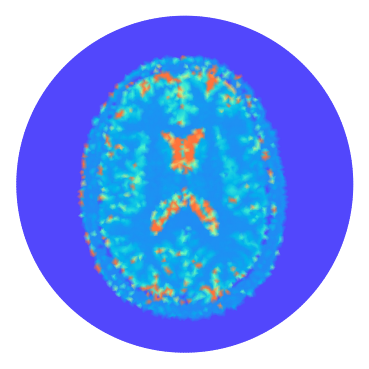}
\caption*{(b)}
\end{subfigure}
\begin{subfigure}{0.05\textwidth}
		\includegraphics[width=1.1cm,height = 3.8cm]{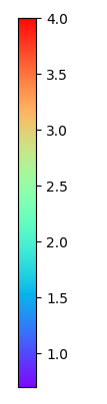}
		\caption*{ }
	\end{subfigure}
\caption{The exact conductivity $\sigma^\dagger$ for different phantoms. (a) the geometrical phantom. (b) actual medical image interpolate into the grid node of the finite element.}
\label{fig:ex2_exact}
\end{figure}

\begin{figure}[H] 
\centering
  \begin{subfigure}{0.22\textwidth}
	\includegraphics[width=3.2cm]{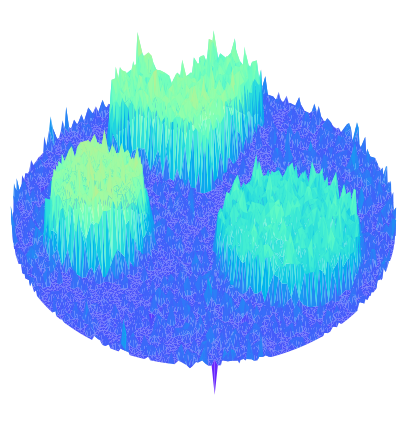}
	\caption*{(a)}
	\end{subfigure}
   \begin{subfigure}{0.22\textwidth}
	\includegraphics[width=3.2cm]{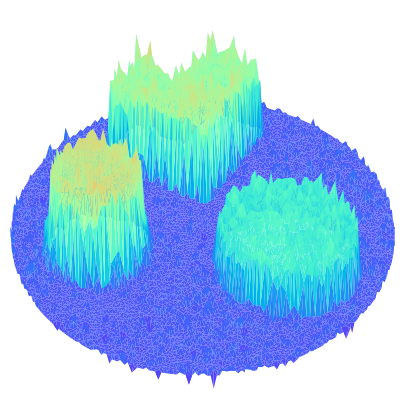}
	\caption*{(b)}
	\end{subfigure}
   \begin{subfigure}{0.22\textwidth}
	\includegraphics[width=3.2cm]{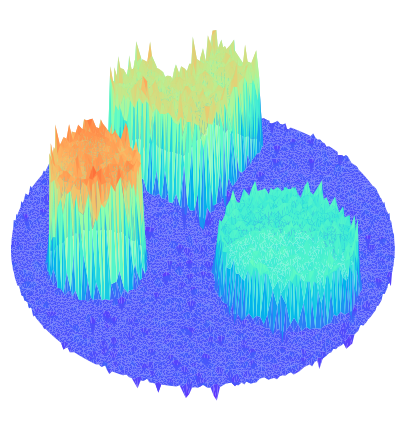}
	\caption*{(c)}
	\end{subfigure}
   \begin{subfigure}{0.22\textwidth}
	\includegraphics[width=3.2cm]{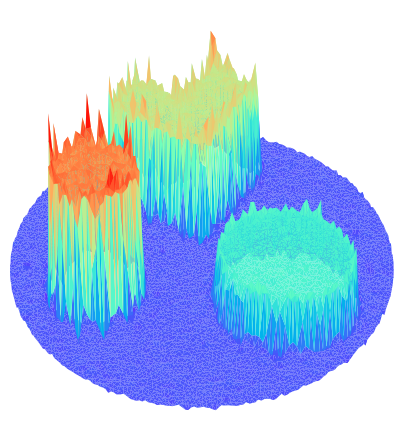}
	\caption*{(d)}
	\end{subfigure}
     \begin{subfigure}{0.08\textwidth}
        \includegraphics[width=0.8cm,height = 3.2cm]{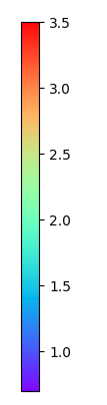} 
        \caption*{ }
    \end{subfigure}

      \begin{subfigure}{0.22\textwidth}
	\includegraphics[width=3.2cm]{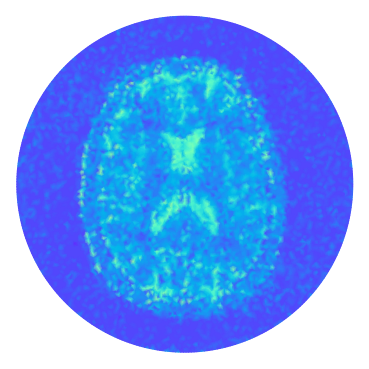}
	\caption*{(e)}
	\end{subfigure}
   \begin{subfigure}{0.22\textwidth}
	\includegraphics[width=3.2cm]{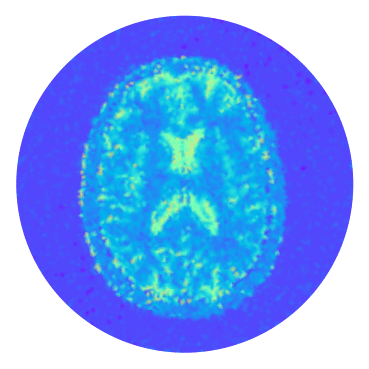}
	\caption*{(f)}
	\end{subfigure}
   \begin{subfigure}{0.22\textwidth}
	\includegraphics[width=3.2cm]{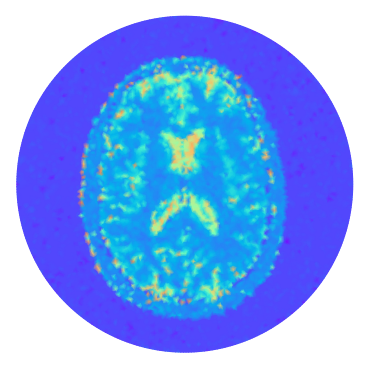}
	\caption*{(g)}
	\end{subfigure}
   \begin{subfigure}{0.22\textwidth}
	\includegraphics[width=3.2cm]{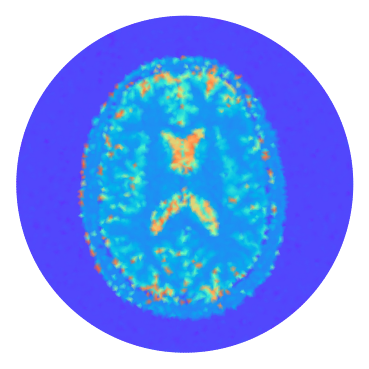}
	\caption*{(h)}
	\end{subfigure}
\begin{subfigure}{0.08\textwidth}
		\includegraphics[width=0.8cm,height = 3.2cm]{brain/ex2_sigma_colorbar.png}
		\caption*{ }
	\end{subfigure}

  \caption{(a)-(h) The reconstructions for both geometrical and brain phantoms under different noise levels under $H_1^\d,H_2^\d,H_3^\d$. (a),(e) $8\%$ noise. (b) (f) $4\%$ noise. (c) (g) $2\%$ noise. (d) (h) $0.8\%$ noise. } 
     \label{fig:ex2_results}
 \end{figure}

We consider two different phantoms. 
\begin{enumerate}
\item  Geometrical shape phantom. The graph has a uniform background value of $1$ and three protruding geometry, one peak corresponding to an elliptical cylinder of height $1.8$,  the other one peak corresponding to a cylinder of height $3$,  another is a concave polyhedron with a peak of height $2.5$, see the Figure \ref{fig:ex2_exact} (a).
\item  Brain phantom\cite{marino2021conductivity}, the accurate conductivity is in the Figure \ref{fig:ex2_exact} (b).  The middle of the image is the linear sulcus, the lateral side is the skull, the area wrapped by the skull is the ventricles, and there are folded meninges in the ventricles.
\end{enumerate}
We may choose the $\Theta_{L^1}$ to execute the reconstruction. The initial subgradient in Algorithm \ref{A1} are chosen by $\xi_{-1} = \xi_0 = 1$, the combination parameter $\lambda_{n,i}^\delta$ are determined by \eqref{eq:choice} with a constant $\alpha = 3$, the step size $\mu_{n,i}^\delta$ is determined by \eqref{eq:mu}, with $\mu_0 = 1.8(1-1/\tau)$ and $\tau = 1.05$. The reconstruction results under different noise levels are plotted in Figure \ref{fig:ex2_results}. It is obvious that, the $L^1$ reconstruction shows its superiority in identifying interfaces and heights of different geometrical shapes. At higher noise levels, the peak surface has more hair thorns, however, the errors on gradient vertical plane and peak horizontal plane are not large. At lower noise levels, the image has less burrs. The same phenomenon also happens in brain phantom. 
 
 In order to show the effectiveness of the method, we compare the TPG-$\Theta$ algorithm in Kaczmarz type with classical Landweber iteration of Kaczmarz type \cite{jin2013landweber}, in which the combination parameter $\lambda_{n,i}^\delta$ are always be zero.  To illustrate the quality of the reconstruction results, we also record the relative reconstruction errors
\begin{align*}
  e_{L^{1}} = \frac{\| \sigma_{n_{\delta}}^{\delta} - \sigma^{\dagger} \Vert_{L^{1}}}{\|\sigma^{\dagger }\Vert_{L^{1}} }\,,\quad
        e_{TV} = \frac{| |\sigma_{n_{\delta}}^{\delta}|_{TV} - |\sigma^{\dagger}|_{TV} |} {|\sigma^{\dagger }|_{TV} } \,,
\end{align*}
and the value of PSNR (peak signal-to-noise ratio). Here the PSNR is defined by
\begin{align*}
\textrm{PSNR} = 10\cdot \log_{10}\frac{\textrm{MAX}^2}{\textrm{MSE}}(\textrm{dB})\,,
\end{align*}
in which the MSE stands for the mean-squared-error per pixel and MAX is the largest value of all pixels. The comparison in Table \ref{table:ex2_error} illustrates that, under the same noise level and stop criteria, the iterative numbers of  TPG-$\Theta$ for Kaczmarz type are significantly reduced, and the relative errors are smaller.

\begin{figure}[H] 
\centering
  \begin{subfigure}{0.49\textwidth}
	\includegraphics[width=7.7cm]{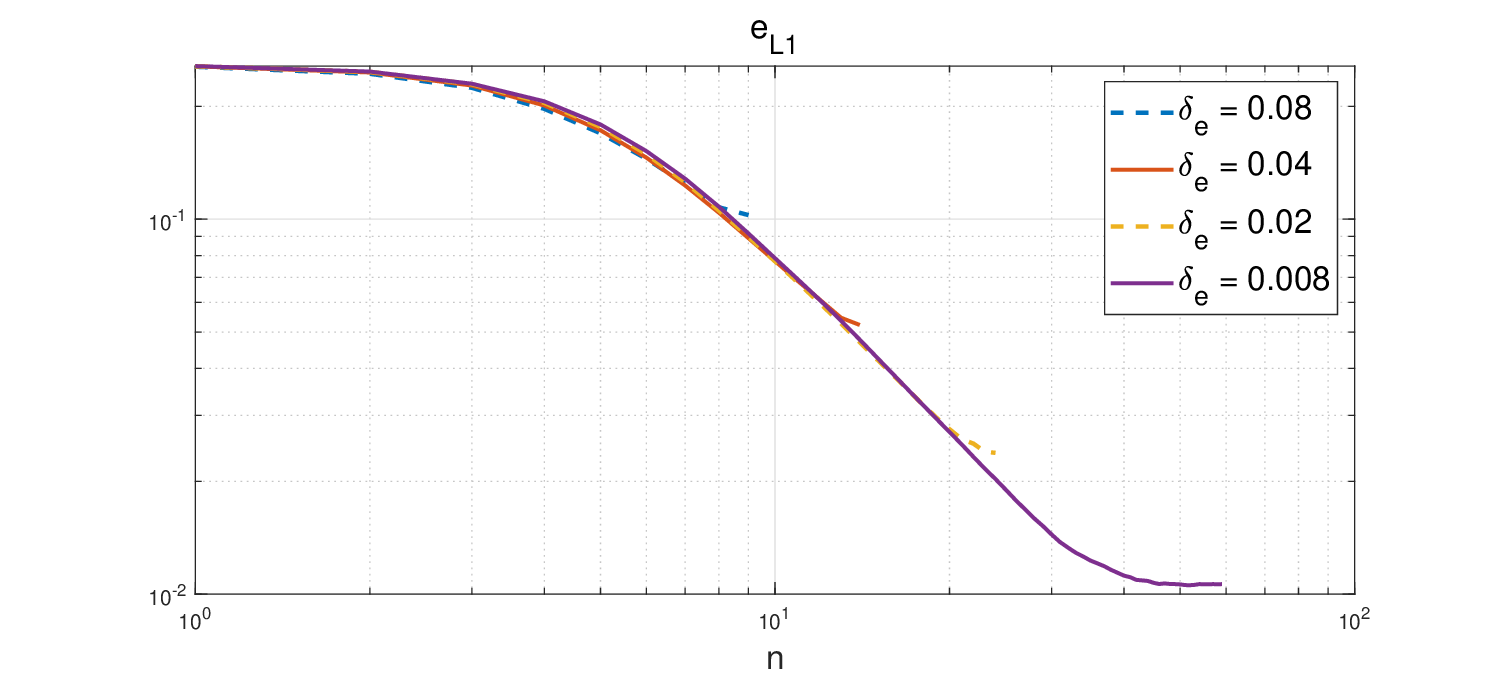}
	\caption*{(a)}
	\end{subfigure}
	\begin{subfigure}{0.49\textwidth}
		\includegraphics[width=7.7cm]{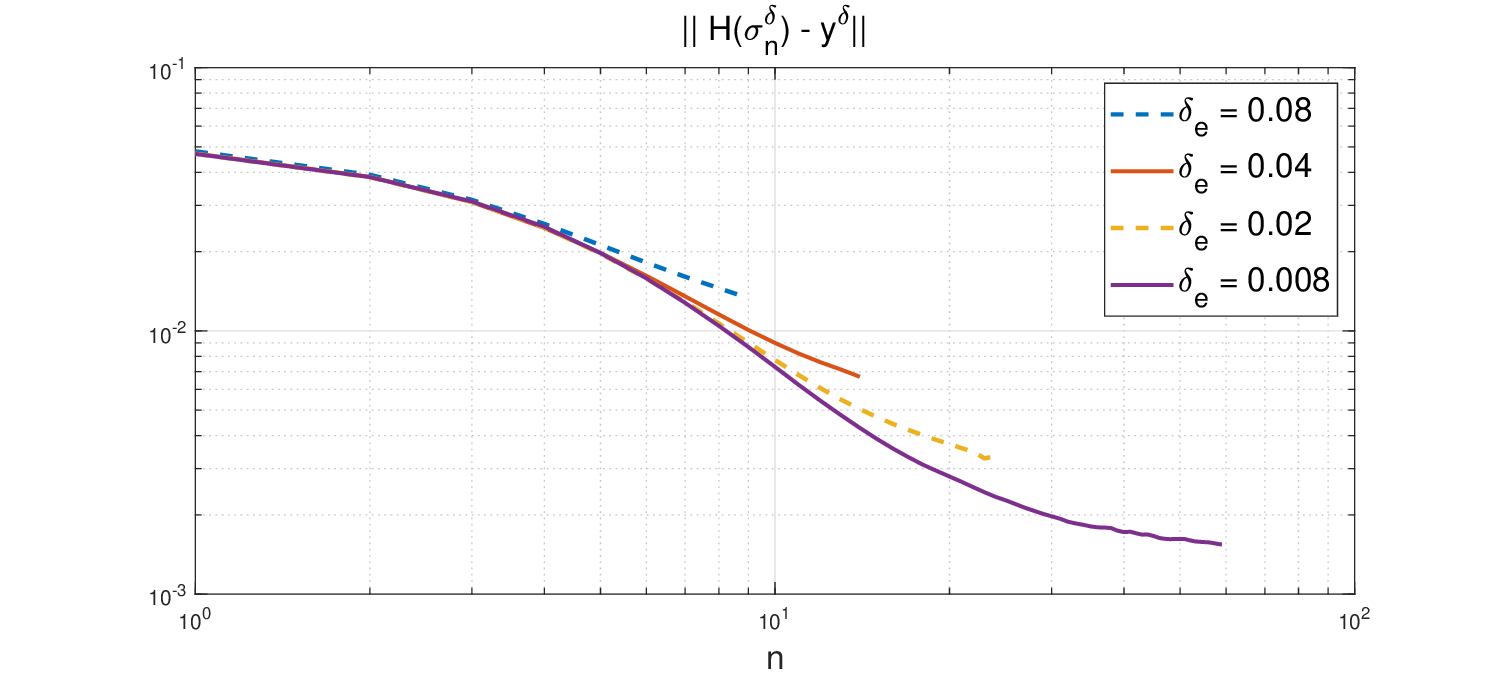}
		\caption*{(b)}
	\end{subfigure}

   \begin{subfigure}{0.49\textwidth}
	\includegraphics[width=7.7cm]{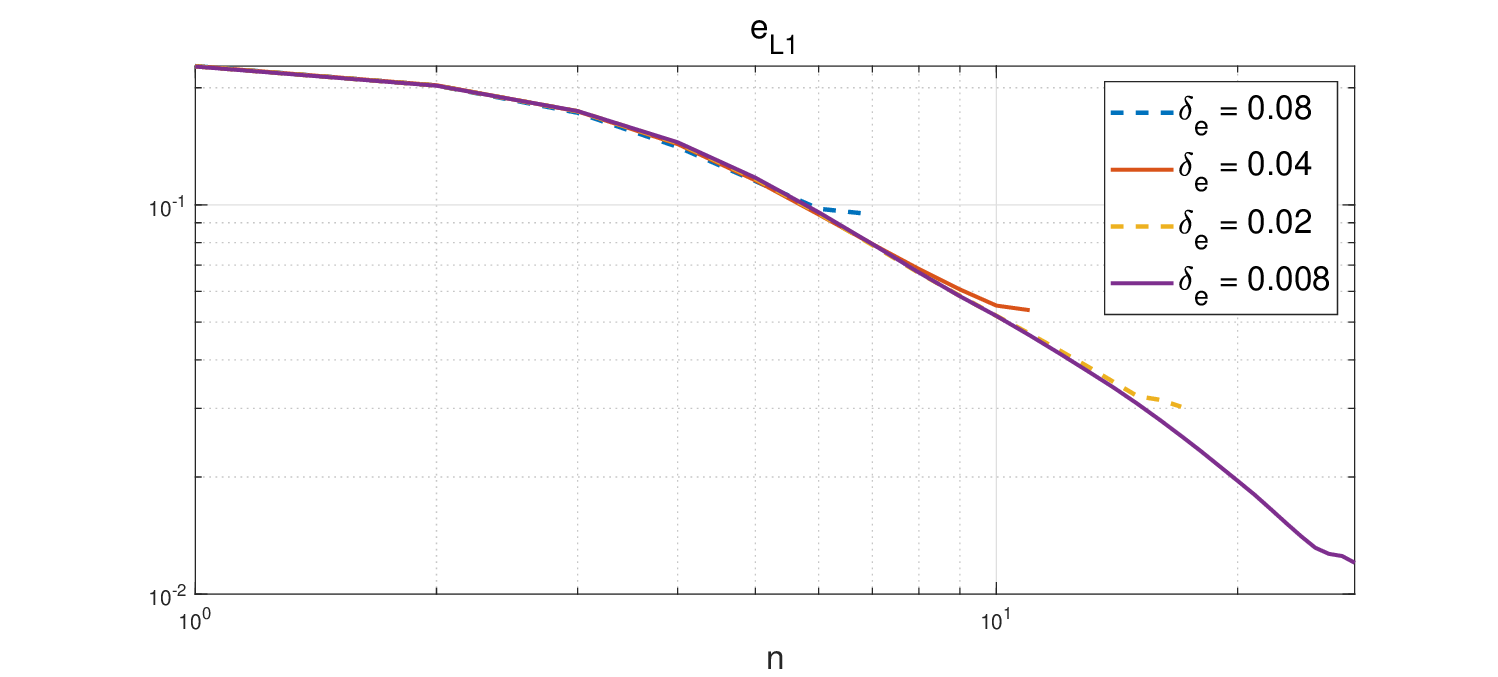}
	\caption*{(c)}
	\end{subfigure}
	\begin{subfigure}{0.49\textwidth}
		\includegraphics[width=7.7cm]{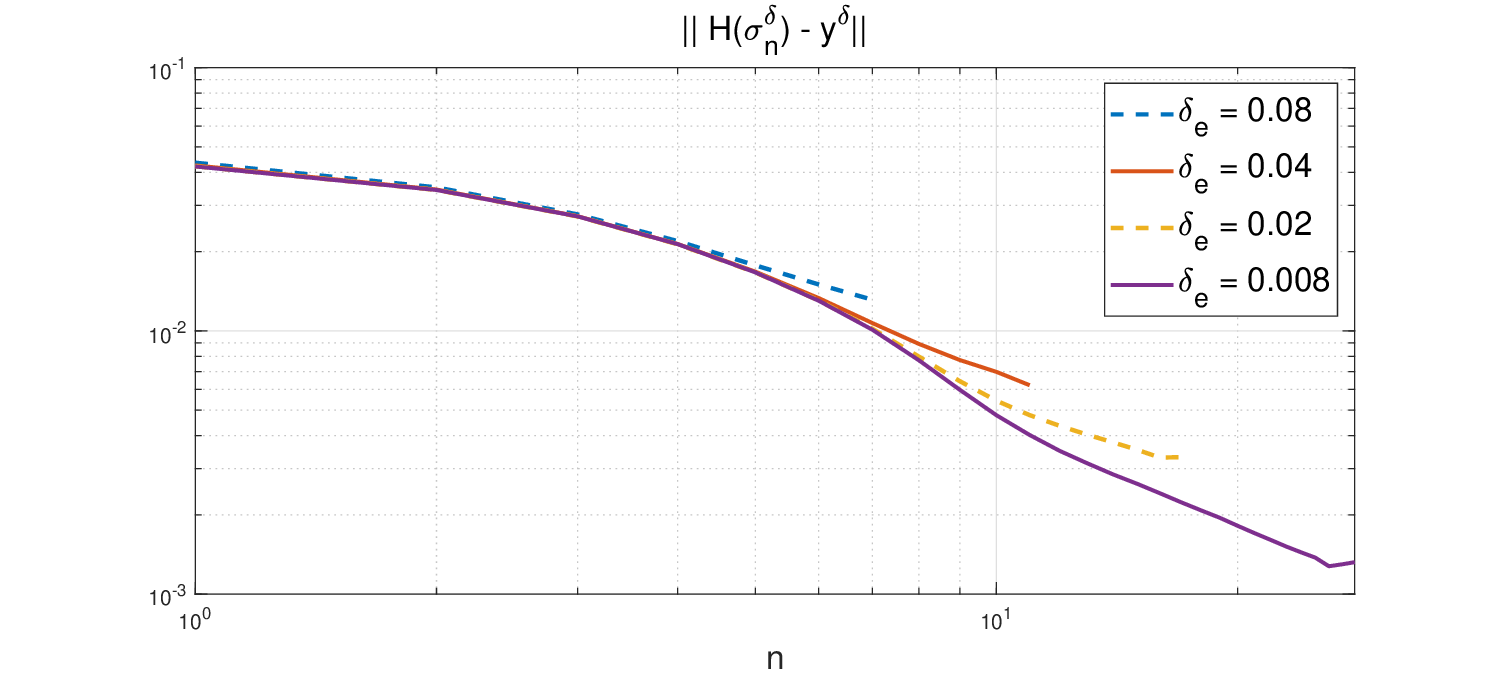}
		\caption*{(d)}
	\end{subfigure}
   \caption{The TPG-$\Theta_{L^1}$-Kaczmarz algorithm reconstruction errors and the residual $\|H(\sigma_n^\d)-y^\d\|$ for the different phantoms. (a)-(b) the geometrical phantom. (c)-(d) the brain phantom.} 
  \label{fig:ex2_L1_curve}
 \end{figure}

\begin{table}[H]
    \centering
    \caption{The comparison for the reconstructions of different phantoms between TPG-$\Theta_{L^1}$-Kaczmarz and Landweber iteration.}
    \label{table:ex2_error}
\begin{tabular}{p{0.12\textwidth}|p{0.07\textwidth}p{0.25\textwidth}p{0.07\textwidth}p{0.1\textwidth}p{0.1\textwidth}p{0.1\textwidth}}
        \specialrule{.1em}{.05em}{.05em} 
            Phantom&     $\delta_e$  & Methods & $n_{\delta}$  & $e_{L^{1}}$ & $e_{TV}$ & PSNR \\

        \specialrule{.1em}{.05em}{.05em} 
        \multirow{9}{*}{Geometry}  &\multirow{2}{*}{$8\%$} &Landweber &18   & 0.113537&0.707887 & 19.9576\\
              &              &TPG-$\Theta_{L^1}$-Kacamarz  &9 &0.102547&0.814492  & 20.6702  \\
        \cline{2-7}
             &   \multirow{2}{*}{$4\%$} &Landweber &40  &0.058384&0.295379 & 24.3346 \\
               &             &TPG-$\Theta_{L^1}$-Kacamarz &14 & 0.052238&0.377399& 25.0684 \\
        \cline{2-7}        
             &\multirow{2}{*}{$2\%$} &Landweber &93 &0.027631&0.093600& 29.3239 \\
                &            &TPG-$\Theta_{L^1}$-Kacamarz &24 &0.023816&0.208704& 30.3543\\
               \cline{2-7}
              &\multirow{2}{*}{$0.8\%$} &Landweber &325 &0.010945&0.027970& 34.6457 \\
                 &           &TPG-$\Theta_{L^1}$-Kacamarz &59 &0.010639&0.122676& 35.6508 \\
        \hline
        \multirow{9}{*}{Brain}  &\multirow{2}{*}{$8\%$} &Landweber &12  & 0.100470&0.166773& 21.0317 \\
        &                          &TPG-$\Theta_{L^1}$-Kacamarz &7 &0.094570&0.068491& 21.5578\\
        \cline{2-7}
                 &\multirow{2}{*}{$4\%$} &Landweber &25 &0.057492&0.137615 & 24.4536\\
        &                          &TPG-$\Theta_{L^1}$-Kacamarz  &11&0.053648&0.084170& 24.9965 \\
        \cline{2-7}
                    &\multirow{2}{*}{$2\%$} &Landweber &51  &0.032920&0.088393 & 28.2113\\
        &                          &TPG-$\Theta_{L^1}$-Kacamarz  &17 &0.030312&0.021722& 29.2629\\
        \cline{2-7}
                     &\multirow{2}{*}{$0.8\%$} &Landweber &120 &0.014276&0.030616  & 33.7022\\
        &                          &TPG-$\Theta_{L^1}$-Kacamarz  &28&0.012044& 0.010117 & 35.5322\\
        \specialrule{.1em}{.05em}{.05em} 
    \end{tabular}
\end{table}
We plot the residual curves $\|H(\sigma_n^\d)-y^\d\|$ and the $L^1$ reconstruction errors with respect to $n$, see the Figure \ref{fig:ex2_L1_curve}. Prior to satisfy the discrepancy principle, the residual $\|H(\sigma_n^\d)-y^\d\|$ and reconstruction errors $e_{L^1}$ decrease while the iterative number $n$ increases.
 
\begin{figure}[H] 
\centering
  \begin{subfigure}{0.22\textwidth}
	\includegraphics[width=3.2cm]{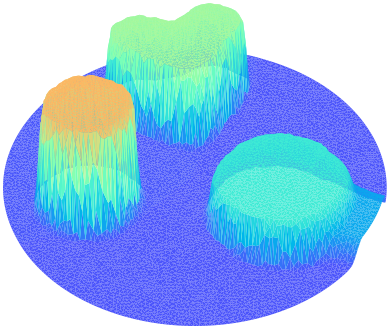}
	\caption*{(a)}
	\end{subfigure}
   \begin{subfigure}{0.22\textwidth}
	\includegraphics[width=3.2cm]{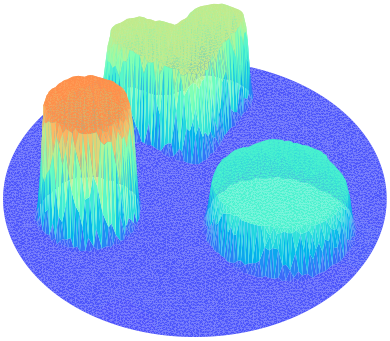}
	\caption*{(b)}
	\end{subfigure}
   \begin{subfigure}{0.22\textwidth}
	\includegraphics[width=3.2cm]{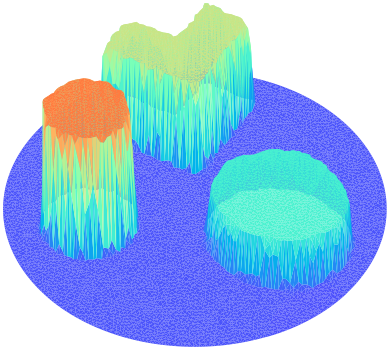}
	\caption*{(c)}
	\end{subfigure}
   \begin{subfigure}{0.22\textwidth}
	\includegraphics[width=3.2cm]{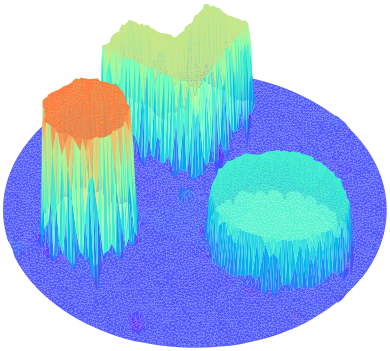}
	\caption*{(d)}
	\end{subfigure}
     \begin{subfigure}{0.08\textwidth}
        \includegraphics[width=0.8cm,height = 3.2cm]{geo/ex2_sig_colorbar.png} 
        \caption*{ }
    \end{subfigure}

      \begin{subfigure}{0.22\textwidth}
	\includegraphics[width=3.2cm]{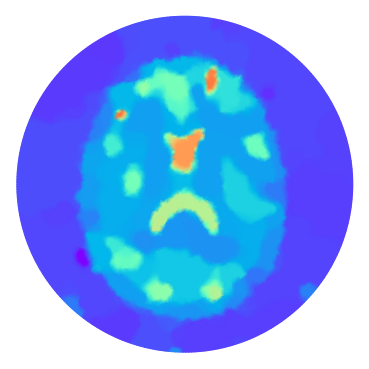}
	\caption*{(e)}
	\end{subfigure}
   \begin{subfigure}{0.22\textwidth}
	\includegraphics[width=3.2cm]{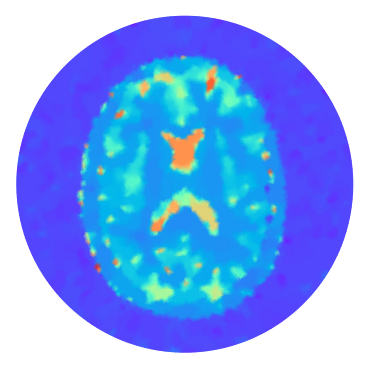}
	\caption*{(f)}
	\end{subfigure}
   \begin{subfigure}{0.22\textwidth}
	\includegraphics[width=3.2cm]{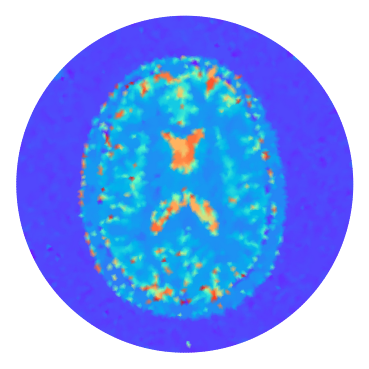}
	\caption*{(g)}
	\end{subfigure}
   \begin{subfigure}{0.22\textwidth}
	\includegraphics[width=3.2cm]{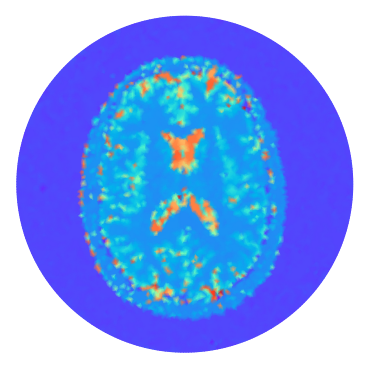}
	\caption*{(h)}
	\end{subfigure}
\begin{subfigure}{0.08\textwidth}
		\includegraphics[width=0.8cm,height = 3.2cm]{brain/ex2_sigma_colorbar.png}
		\caption*{ }
	\end{subfigure}

  \caption{(a)-(h) The reconstructions for both geometrical and brain phantoms under different levels of noise under $H_{1-3}^\d$. (a),(e) $8\%$ noise. (b) (f) $4\%$ noise. (c) (g) $2\%$ noise. (d) (h) $0.8\%$ noise.  } 
  \label{fig:ex2_TV}
 \end{figure}

We also test the $\Theta_{TV}$ to execute the algorithms. The results are plotted in Figure \ref{fig:ex2_TV} and the details are recorded in Table \ref{tab:ex2_TV}. By contrasting reconstruction images and errors with previous $\Theta_{L^1}$ penalty, the results are more satisfactory. In addition, compared with classical Landweber in Kaczmarz type, the TPG-$\Theta_{TV}$-Kaczmarz algorithm significantly reduces the iteration number, and save the storage. However, due the high computational costs in calculating the minimization problem \eqref{eq:ROF}, the calculation consumption has increased. The accuracy losses in solving the minimization problem are also reflected in the final results.

\begin{figure}[H] 
\centering
  \begin{subfigure}{0.49\textwidth}
	\includegraphics[width=7.7cm]{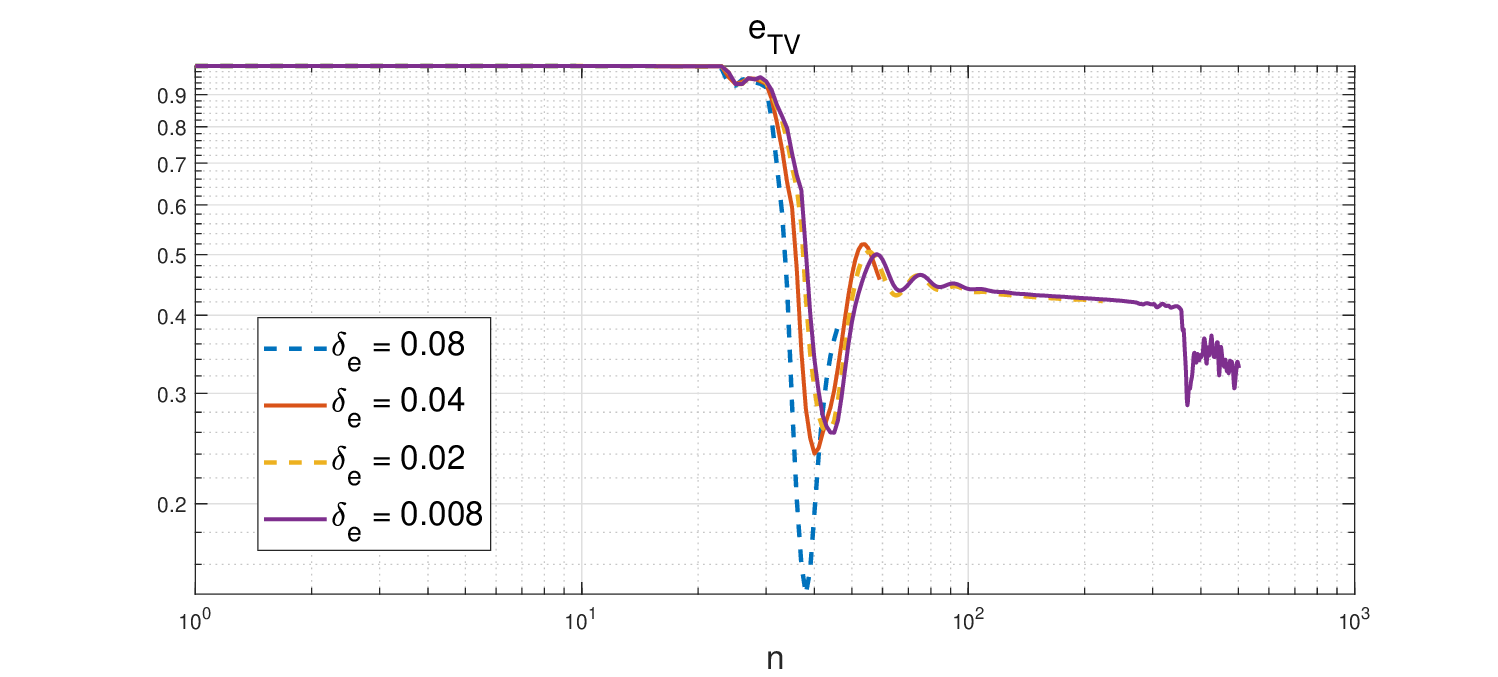}
	\caption*{(a)}
	\end{subfigure}
	\begin{subfigure}{0.49\textwidth}
		\includegraphics[width=7.7cm]{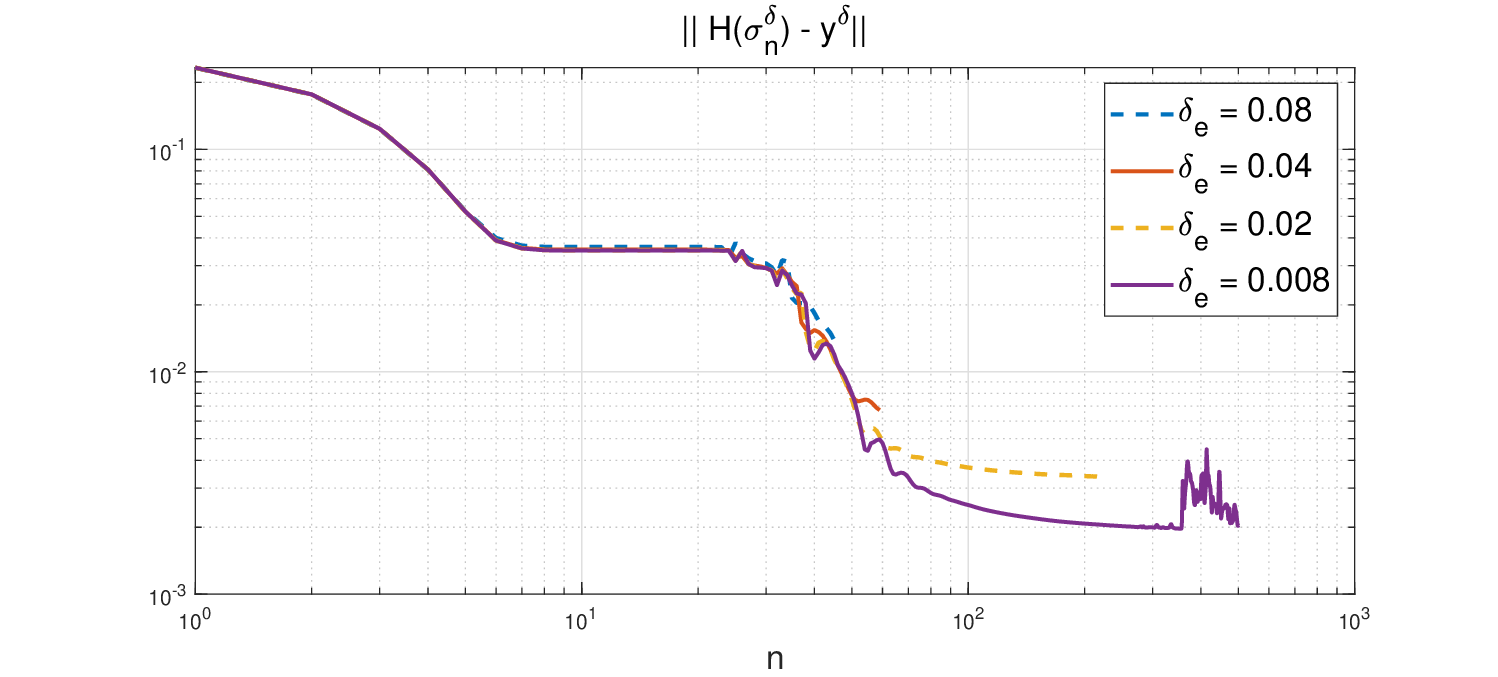}
		\caption*{(b)}
	\end{subfigure}

   \begin{subfigure}{0.49\textwidth}
	\includegraphics[width=7.7cm]{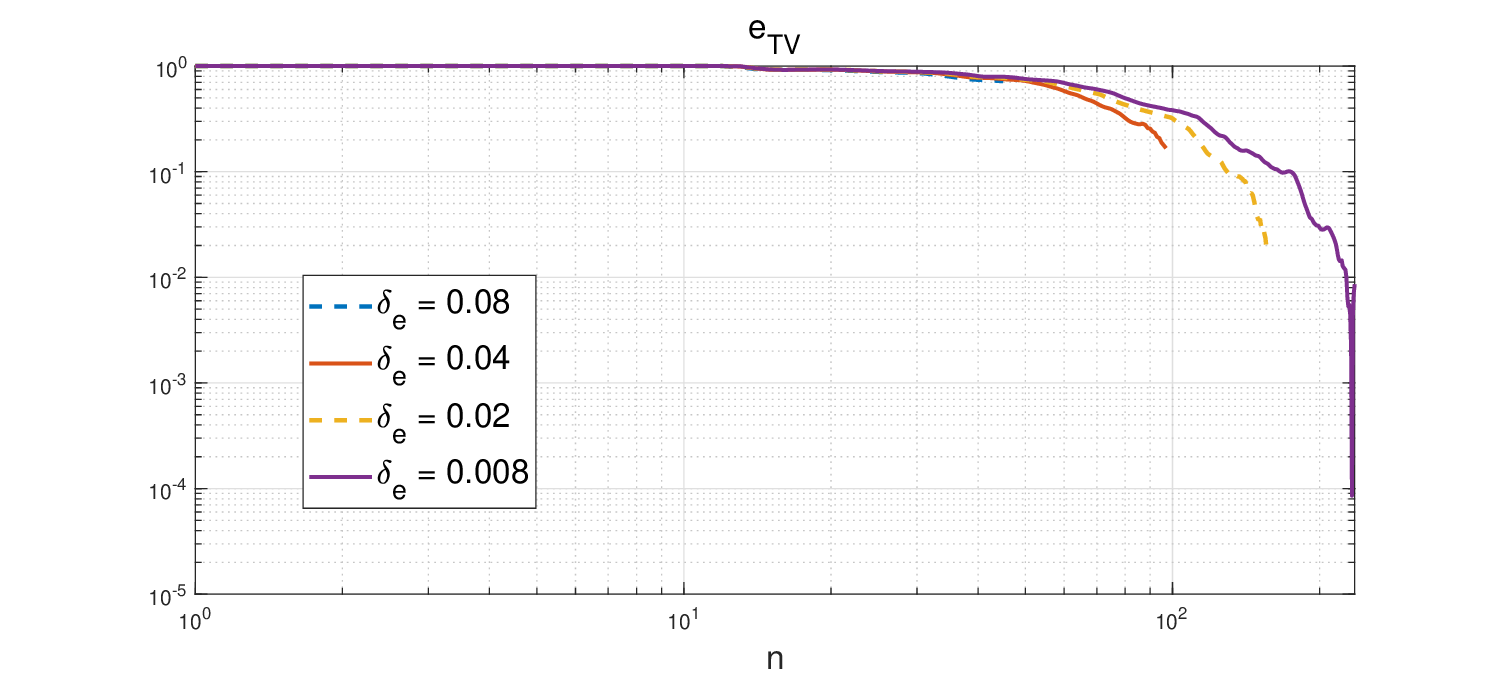}
	\caption*{(c)}
	\end{subfigure}
	\begin{subfigure}{0.49\textwidth}
		\includegraphics[width=7.7cm]{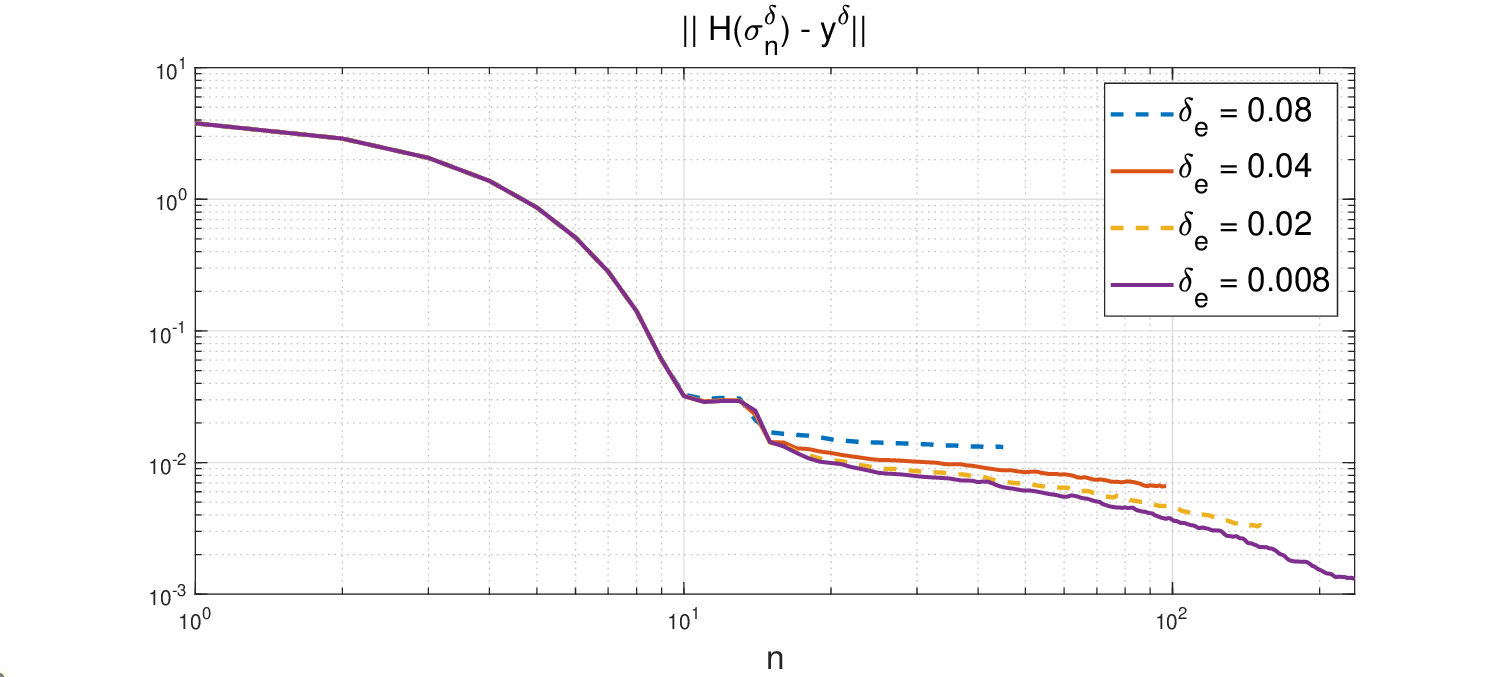}
		\caption*{(d)}
	\end{subfigure}
   \caption{The TPG-$\Theta_{TV}$-Kaczmarz algorithm reconstruction errors and the residual $\|H(\sigma_n^\d)-y^\d\|$ for the different phantoms. (a)-(b) the geometrical phantom. (c)-(d) the brain phantom.} 
  \label{fig:ex2_TV_curve}
 \end{figure}

 Figure \ref{fig:ex2_TV_curve} plots the residual curves $\|H(\sigma_n^\d)-y^\d\|$ and the $TV$ reconstruction errors with respect to $n$. There are some small jumps in the $e_{TV} $ during the iteration, which related to the difference in total variation. This seems that the iterative calculation locally processes discontinuities, making the discontinuities and discontinuous textures more clearer. Thus it is mostly be the tail of the plot, where the $e_{TV}$ decreases, that is of significance. The plot shows obviously that the convergence and the stability of the algorithm even in the noisy case.

\begin{table}[H]
    \centering
    \caption{The comparison for the reconstructions of different phantoms between TPG-$\Theta_{TV}$-Kaczmarz and Landweber iteration.}
    \label{tab:ex2_TV}
\begin{tabular}{p{0.12\textwidth}|p{0.07\textwidth}p{0.25\textwidth}p{0.07\textwidth}p{0.1\textwidth}p{0.1\textwidth}p{0.1\textwidth}}
        \specialrule{.1em}{.05em}{.05em} 
            Phantom&     $\delta_e$  & Methods & $n_{\delta}$  & $e_{L^{1}}$ & $e_{TV}$ & PSNR \\

        \specialrule{.1em}{.05em}{.05em} 
        \multirow{9}{*}{Geometry}  &\multirow{2}{*}{$8\%$} &Landweber &254   & 0.055821&0.528660 & 23.5124\\
              &              &TPG-$\Theta_{TV}$-Kacamarz  &46 &0.056442&0.383817 &  23.4731  \\
        \cline{2-7}
             &   \multirow{2}{*}{$4\%$} &Landweber &617  &0.022002&0.465395 & 27.0950\\
               &             &TPG-$\Theta_{TV}$-Kacamarz &59 & 0.027035&0.456140& 27.1813 \\
        \cline{2-7}        
             &\multirow{2}{*}{$2\%$} &Landweber &8149&0.013069&0.422790& 31.2199 \\
                &            &TPG-$\Theta_{TV}$-Kacamarz &223 &0.014281&0.421279& 31.1811\\
               \cline{2-7}
              &\multirow{2}{*}{$0.8\%$} &Landweber &16323 &0.013534 & 0.415861 &31.5072 \\
                 &           &TPG-$\Theta_{TV}$-Kacamarz &480 &0.014456&0.337332 &33.9794 \\
        \hline
        \multirow{9}{*}{Brain}  &\multirow{2}{*}{$8\%$} &Landweber &276  & 0.077656&0.797530& 21.0617 \\
        &                          &TPG-$\Theta_{TV}$-Kacamarz &45&0.080472&0.719443& 21.2758\\
        \cline{2-7}
                 &\multirow{2}{*}{$4\%$} &Landweber &960 &0.044854&0.487634 & 24.2998\\
        &                          &TPG-$\Theta_{TV}$-Kacamarz  &97&0.047308&0.167007& 23.8465 \\
        \cline{2-7}
                    &\multirow{2}{*}{$2\%$} &Landweber &2281  &0.024636&0.260189 & 27.8652\\
        &                          &TPG-$\Theta_{TV}$-Kacamarz  &157 &0.025804&0.017585& 29.3383\\
        \cline{2-7}
                     &\multirow{2}{*}{$0.8\%$} &Landweber &5910 &0.011146&0.103770  & 32.2589\\
        &                          &TPG-$\Theta_{TV}$-Kacamarz  &236&0.011404& 0.008582 & 31.6143 \\
        \specialrule{.1em}{.05em}{.05em} 

    \end{tabular}
\end{table}

\subsection{Numerical reconstructions for partial interior data}
\begin{figure}[!htbp]    


        \begin{subfigure}{0.22\textwidth}
      \includegraphics[width=3.2cm]{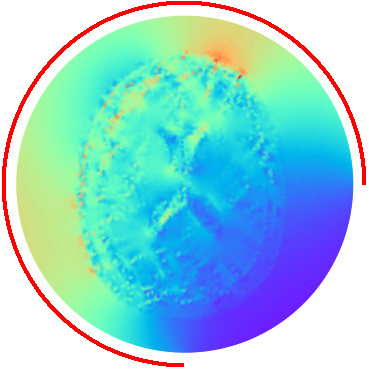}
      \end{subfigure}
      \begin{subfigure}{0.22\textwidth}
              \includegraphics[width=3.2cm]{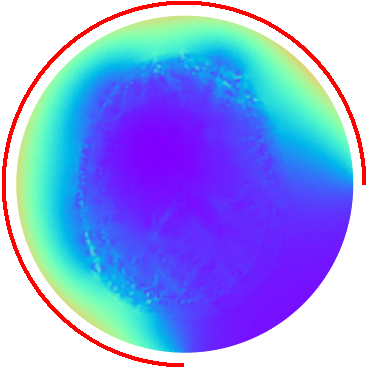}
      \end{subfigure}
      \begin{subfigure}{0.22\textwidth}
              \includegraphics[width=3.2cm]{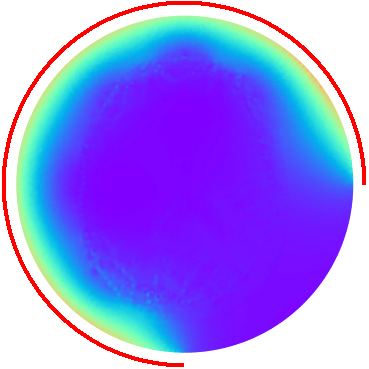}
      \end{subfigure}
      \begin{subfigure}{0.22\textwidth}
          \includegraphics[width=3.2cm]{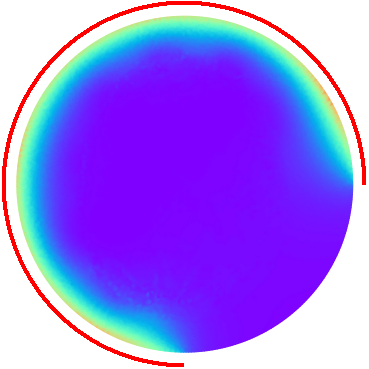}
      \end{subfigure}
      \begin{subfigure}{0.08\textwidth}
              \includegraphics[width=0.5cm,height = 3.2cm]{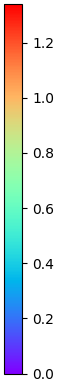}
      \end{subfigure}

           \begin{subfigure}{0.22\textwidth}
             \includegraphics[width=3.2cm]{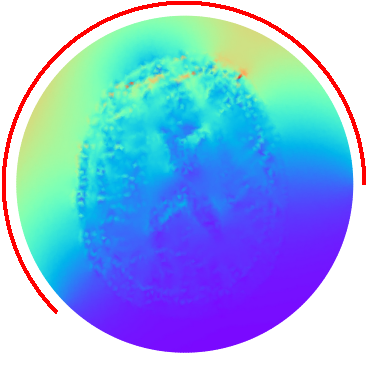}
       \end{subfigure}
       \begin{subfigure}{0.22\textwidth}
             \includegraphics[width=3.2cm]{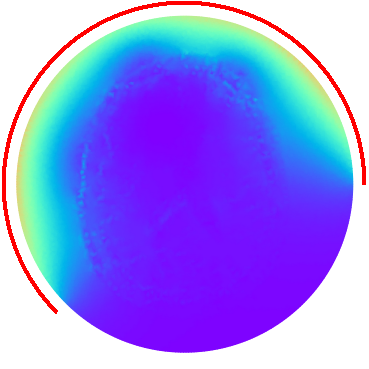}
       \end{subfigure}
       \begin{subfigure}{0.22\textwidth}
             \includegraphics[width=3.2cm]{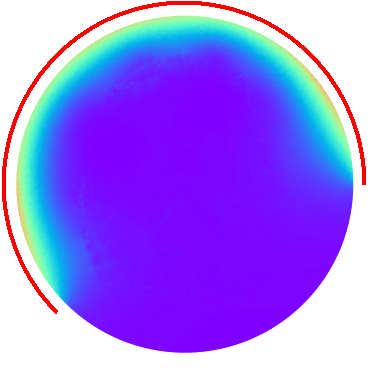}
       \end{subfigure}
       \begin{subfigure}{0.22\textwidth}
         \includegraphics[width=3.2cm]{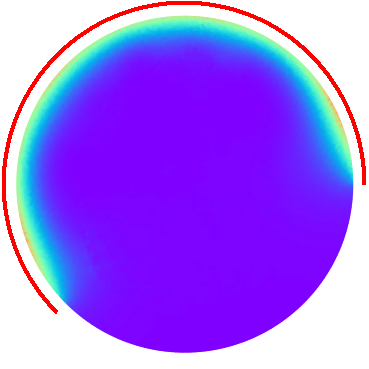}
       \end{subfigure}
       \begin{subfigure}{0.08\textwidth}
             \includegraphics[width=0.5cm,height = 3.2cm]{limited_angle/obs_colorbar.png}
       \end{subfigure}

       \begin{subfigure}{0.22\textwidth}
             \includegraphics[width=3.2cm]{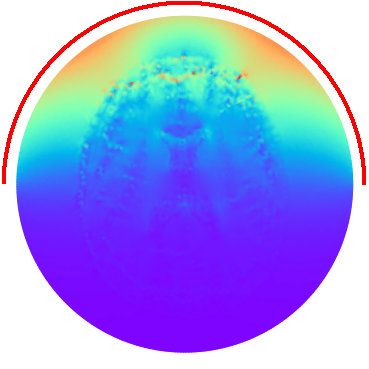}
       \end{subfigure}
       \begin{subfigure}{0.22\textwidth}
             \includegraphics[width=3.2cm]{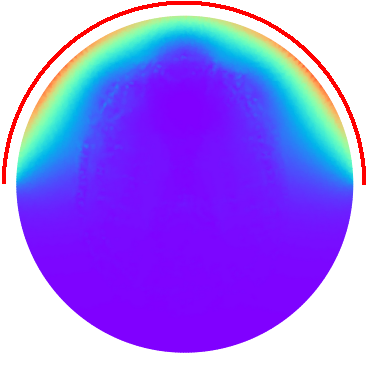}
       \end{subfigure}
       \begin{subfigure}{0.22\textwidth}
             \includegraphics[width=3.2cm]{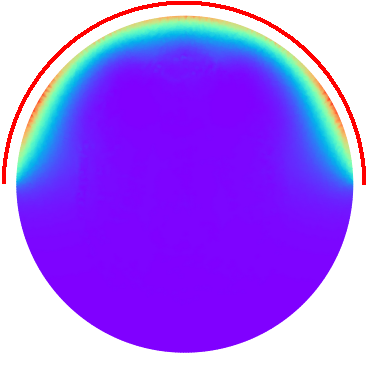}
       \end{subfigure}
       \begin{subfigure}{0.22\textwidth}
         \includegraphics[width=3.2cm]{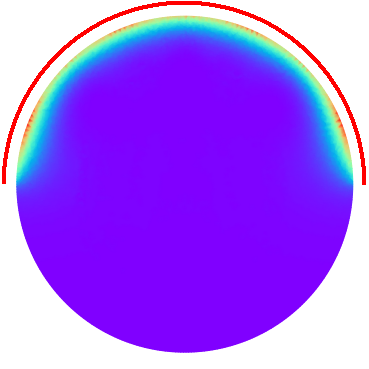}
       \end{subfigure}
       \begin{subfigure}{0.08\textwidth}
             \includegraphics[width=0.5cm,height = 3.2cm]{limited_angle/obs_colorbar.png}
       \end{subfigure} 

        \caption{Power densities $H_{1-4}$ with boundary data $f_{1-4}(r,\theta)$.  The red curves indicate the available boundary.}
        \label{fig:limited_angle_interior}
\end{figure}

The TPG-$\Theta$ with Kaczmarz type algorithm is able to extend straightforwardly to partial interior data. By modifying the boundary $\Gamma(\alpha)\subset\partial\Omega$ by 
$$\Gamma(\alpha) = (r,\theta)\in \frac{1}{2}\times [0,\alpha]\,,$$
the $\partial\Omega\setminus\Gamma$ represents the occluded area, that is, boundary information and power densities are not available. We set
\begin{align}\label{eq:circular}
f_i(r,\theta) = \sin\left(\frac{2i\pi\theta}{\alpha}\right)\,,\quad\forall(r,\theta)\in\Gamma(\alpha),\quad i=1,2,3,4\,
\end{align}
and always asssume $f_i=0$ on the remaining occluded part.  The power densities are depicted in Figure \ref{fig:limited_angle_interior} for the different angles. The circle (segment) in the figures indicates the available boundary. The internal structures are not clearly, especially in the high-frequency case. 

\begin{figure}[!htbp]

        \begin{subfigure}{0.22\textwidth}
      \includegraphics[width=3.2cm]{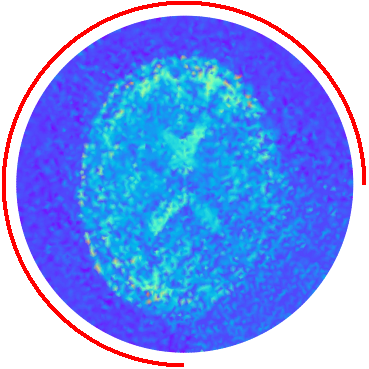}
    \caption*{(a) $\delta_e = 8\%$.}
      \end{subfigure}
        \begin{subfigure}{0.22\textwidth}
      \includegraphics[width=3.2cm]{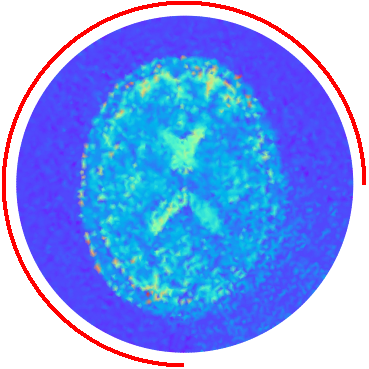}
    \caption*{(b) $\delta_e = 4\%$.}
      \end{subfigure}
             \begin{subfigure}{0.22\textwidth}
      \includegraphics[width=3.2cm]{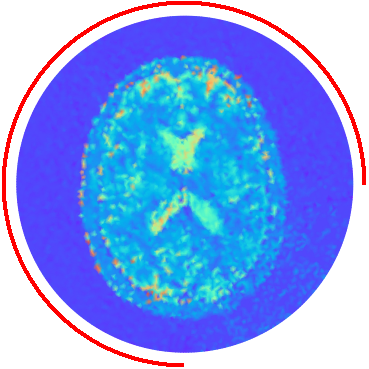}
    \caption*{(c) $\delta_e = 2\%$.}
      \end{subfigure}
              \begin{subfigure}{0.22\textwidth}
      \includegraphics[width=3.2cm]{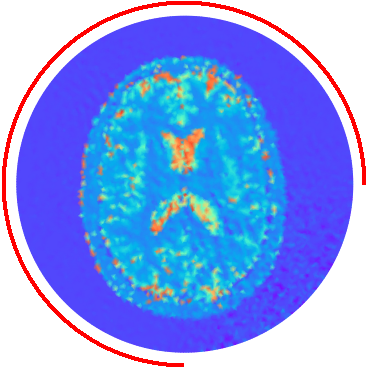}
    \caption*{(d)$\delta_e = 0.8\%$.} 
      \end{subfigure}
      \begin{subfigure}{0.08\textwidth}
		\includegraphics[width=0.8cm,height = 3.2cm]{brain/ex2_sigma_colorbar.png}
		\caption*{ }
	\end{subfigure}

            \begin{subfigure}{0.22\textwidth}
      \includegraphics[width=3.2cm]{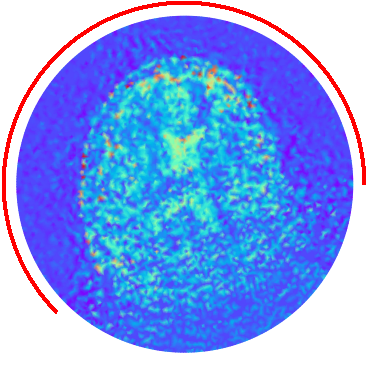}
    \caption*{(a) $\delta_e = 8\%$.}
      \end{subfigure}
        \begin{subfigure}{0.22\textwidth}
      \includegraphics[width=3.2cm]{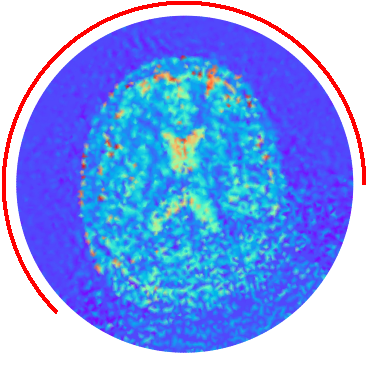}
    \caption*{(b) $\delta_e = 4\%$.}
      \end{subfigure}
             \begin{subfigure}{0.22\textwidth}
      \includegraphics[width=3.2cm]{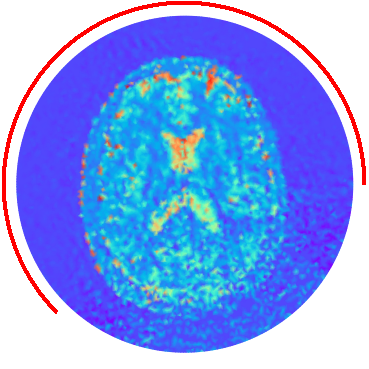}
    \caption*{(c) $\delta_e = 2\%$.}
      \end{subfigure}
              \begin{subfigure}{0.22\textwidth}
      \includegraphics[width=3.2cm]{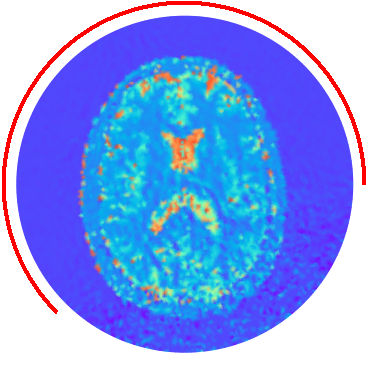}
    \caption*{(d)$\delta_e = 0.8\%$.} 
      \end{subfigure}
      \begin{subfigure}{0.08\textwidth}
		\includegraphics[width=0.8cm,height = 3.2cm]{brain/ex2_sigma_colorbar.png}
		\caption*{ }
	\end{subfigure}

    \begin{subfigure}{0.22\textwidth}
      \includegraphics[width=3.2cm]{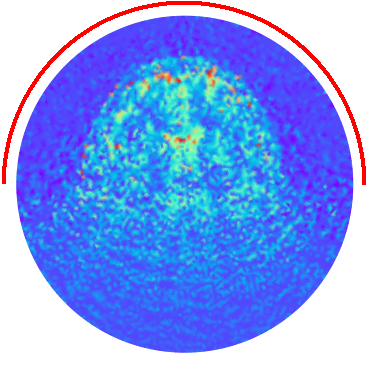}
    \caption*{(e) $\delta_e = 8\%$.}
      \end{subfigure}
       \begin{subfigure}{0.22\textwidth}
      \includegraphics[width=3.2cm]{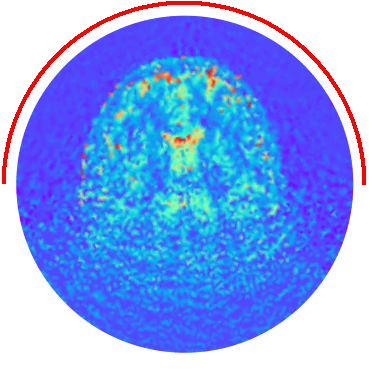}
    \caption*{(f) $\delta_e = 4\%$.}
      \end{subfigure}
             \begin{subfigure}{0.22\textwidth}
      \includegraphics[width=3.2cm]{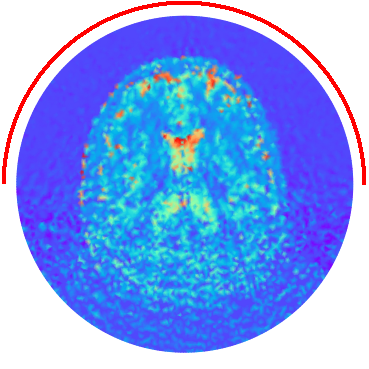}
    \caption*{(g) $\delta_e = 2\%$.}
      \end{subfigure}
              \begin{subfigure}{0.22\textwidth}
      \includegraphics[width=3.2cm]{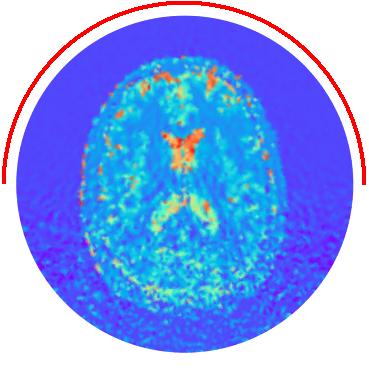}
    \caption*{(h) $\delta_e = 0.8\%$.}
      \end{subfigure}
      \begin{subfigure}{0.08\textwidth}
		\includegraphics[width=0.8cm,height = 3.2cm]{brain/ex2_sigma_colorbar.png}
		\caption*{ }
	\end{subfigure}

 \caption{The reconstructions with limited angles under different noise levels $\delta_e = 8\%,4\%,2\%$ and $ 0.8\%$. The red curves indicate the available boundary. }
       \label{fig:limited_angle_results}
\end{figure}

The reconstructions with different angles and different noise levels $\delta_e$ are plotted in Figure \ref{fig:limited_angle_results}. The inclusions in the phantom are reconstructed accurately near the observed area, but blurry in unobserved area. However, it seems the algorithms are making great
efforts to compensate for the missing conductivity reconstruction outside the data domain. This especially happens in $\alpha = \pi$ case, when the relative noise level is small ($\delta_e = 0.8\%$), although there is almost no information in the remaining angle range, the reconstruction results still provide a rough outline within the occluded area.




\bibliographystyle{siamplain}
\bibliography{sisc_bib}
\end{document}